\documentclass{amsart}
\usepackage{amsmath}
\usepackage{amssymb}
\usepackage{amsthm}
\usepackage{amscd}
\usepackage{epsf}

\newcounter{AbcT}

\newtheorem {Theorem}    {Theorem}[section]
\newtheorem {Definition} {Definition} [section]

\newtheorem {Lemma}      [Theorem]    {Lemma}

\newtheorem {Claim}      [Theorem]    {Claim}

\newtheorem {Conjecture}[Theorem]     {Conjecture}

\newtheorem {Assumption}[Theorem]     {Assumption}

\newcommand {\Heads}[1]   {\smallskip\pagebreak[1]\noindent{\bf
#1{\hskip 0.2cm} }}
\newcommand {\Head}[1]    {\Heads{#1:}}

\newcommand {\proofs}     {\proof}
\newcommand {\prooft}[1]  {\Head{Proof #1}\nopagebreak[2]}

\newcommand  {\QED}    {\def\qedsymbol{$\blacksquare$}\qed}

\newcommand{\ignore}[1]{}

\newcommand{\sign}{{\mbox{{\rm sign}}}}

\newcommand{\poly}{{\mbox{{\rm poly}}}}
\renewcommand{\top}{{\mbox{{\rm top}}}}
\newcommand{\Maj}{{\mbox{{\rm Maj}}}}
\newcommand{\polylog}{{\mbox{{\rm polylog}}}}

 \def\eps{\epsilon} \def\E{{\bf{E}}}
\def\P{{\bf{P}}}

\def\R{\hbox{I\kern-.2em\hbox{R}}}
\def\A{{\mathcal{A}}}
\def\C{{\mathcal{C}}}

\def\|{\, | \, }
\def\v0{{\bf 0}}
\def\one{{\bf 1}}
\def\0{\hat{0}}
\def\1{\hat{1}}

\def\pa{{\tt path}}
\def\path{{\tt path}}

\def\lam{\lambda}

\def\Trees{{\bf T}}
\def\Mut{{\bf M}}

\def\phi{\varphi}

\def\be{\begin{equation}}
\def\ee{\end{equation}}

\def\one{{\bf 1}}

\begin{document}
\title{Phase transitions in Phylogeny}
\author{Elchanan Mossel} 
\address{Department of Statistics, Evans Hall, University of
  California, Berkeley, California 94720-3860}
\email{mossel@stat.berkeley.edu}
\thanks{Supported by a Miller Fellowship. Most of the research
  reported here was conducted while the author was a PostDoc in theory
  group, Microsoft Research.}


\date{Submitted May 21, 2002, in revised form April 10,2003, accepted 
April 10, 2003}




\begin{abstract}

We apply  
the theory of markov random fields on trees to derive
a phase transition in the number of samples needed 
in order to reconstruct phylogenies.

We consider the Cavender-Farris-Neyman model of evolution on trees, 
where all the inner nodes have degree at least $3$, 
and the net transition on each edge is
bounded by $\eps$. Motivated by a conjecture by M. Steel, 
we show that if $2 (1 - 2 \eps)^2 > 1$,
then for balanced trees, the topology of the underlying tree, 
having $n$ leaves, can be reconstructed
from $O(\log n)$ samples (characters) at the leaves. 
On the other hand, we show that if $2 (1 - 2 \eps)^2 < 1$,
then there exist topologies which require at least $n^{\Omega(1)}$ 
samples for reconstruction. 

Our results are the first rigorous results to establish the role
of phase transitions for markov random fields on trees
as studied in probability, statistical physics and information theory to the
study of phylogenies in mathematical biology.

\end{abstract}
\maketitle

\section{Introduction}

Phylogenetic trees commonly model evolution of species.
In this paper we study reconstruction of phylogenetic trees
from samples (characters) of data at the leaves of the tree, 
and the relationship between this problem
and the theory of markov random fields on trees.

We apply tools from the theory of Ising model on trees 
to derive a phase transition in the number of samples needed in order
to reconstruct the topology of a tree for the Cavender-Farris-Neyman 
model of evolution. 

Cavender, Farris and Neyman \cite{Neyman,Cavender,Farris} 
introduced a model of evolution of binary characters.
In this model, the evolution of characters is governed by the Ising
model on the tree of species. It is assumed that characters
evolve identically and independently.
In statistical physics, the study of the Ising model on trees, 
which dates back to the first half of the 20'th century, focused
mostly on regular trees (named ``Bethe lattice'', 
``homogeneous trees `` or ``Cayley trees'' in statistical physics),
and only more recently on general trees \cite{Lyons,Peres}.
However, the problem of reconstructing a tree from samples of data at
the leaves was not studied in this context.


The threshold for the extremality of the free measure for the Ising
model on the tree
will play a crucial role for phylogenies.
The study of this threshold was initiated in \cite{Sp,Hi}.
The exact threshold for regular trees was found in \cite{BRZ},
see also \cite{EKPS,Io,M:recur}.

\subsection{Definitions}
We begin by defining the evolution process.
For a tree $T=(V,E)$ rooted at $\rho \in V$, we direct all edges away from
the root, so that edge $e$ is written as $e=(v,u)$, where $v$ is
on the unique path connecting $\rho$ to $u$.

Let $\A$ be a finite set representing the values of 
some genetic characteristic. Illustrative examples are 
$\A = \{A,C,G,T\}$, or $\A = \{20 \mbox{ amino acids}\}$.
We will often refer to the elements of $\A$ as colors.

The propagation of the genetic character $\sigma$ from $\rho$ to the nodes of
the tree $T$ is modeled in the following manner. 
The root color is chosen according to some
initial distribution $\pi$, so that $\P[\sigma_{\rho} = i] = \pi_i$.
The mutation along edge $e$ is encoded by a stochastic matrix  $(M^e_{i,j})_{i,j=1}^\ell$.
For edge $e=(v,u)$, it holds that $\P[\sigma_u = j | \sigma_v = i] = M^{(v,u)}_{i,j}$.
Moreover, if $\pa(u,v) = \{u = v_0,\ldots,v_{\ell} = v\}$, is the path from $u$ to $v$ in $T$,
and $\Delta(v) = \{ w : v \notin \pa(\rho,w)\}$,
then it is assumed that $(\sigma_v)$ satisfies the following markov
property (see e.g. \cite{Georgii,Liggett} for the general definition
of markov random field):
\[
\P[\sigma_u = j | \sigma_v = i, (\sigma_w)_{w \in \Delta(u)}] =
M^{(v,u)}_{i,j}.
\]

One of the fundamental problems of mathematical and computational
biology is the {\em reconstruction} of phylogenetic trees. Our model
of evolution is defined on a rooted tree.
However, we will only consider the reconstruction of the un-rooted tree
and ignore the problem of reconstructing the root. Indeed, in many
cases it is impossible to reconstruct the root given the data at the
leaves. For example, if all the $M^e$ are reversible with respect to
the same distribution $\pi$ (which is also the initial distribution), 
then without additional data or
assumptions on the model, it is impossible to distinguish a root.     

For a tree $T=(V,E)$, we call $v \in V$ a {\em leaf}, 
if $v$ has degree $1$ in the graph $(V,E)$.
We write $\partial T$ for the {\em boundary} of the tree, i.e., the
set of all leaves of $T$.
 

Let $T=(V,E)$ be a tree on $n$ leaves. 
Consider the evolution process on $T$, where we consider $T$ as a
tree rooted at $\rho \in V$, where $\rho$ is not one of the leaves of
$T$; let $(M^e)_{e \in E}$ be the collection of mutation matrices and 
$(\pi_i)_{i \in \A}$ the initial distribution at $\rho$.
For a coloring $\sigma$ on the vertices of the tree, denote by
$\sigma_{\partial T}$ the values of the color at the 
boundary of the tree.
Suppose that $k$ independent samples of the above process, 
$(\sigma^t_v)_{1 \leq t \leq k;v \in T}$ are given.
In biology it is common to call these samples {\em characters}, 
we refer to these either as {\em samples} or as characters.
The objective is to find $T$ given the samples at the leaves 
$(\sigma^t_{\partial T})_{t=1}^k$.

The standard assumption in phylogeny is that all the internal 
degrees in $T$ are $3$; it is also assumed that all rooted trees are
rooted at internal vertices, see, e.g., \cite{ESSW1,SSH}. 
We will slightly relax the first assumption.
\begin{Assumption}
We assume that the evolution process is defined on rooted trees $T$, such
that  all internal degrees of $T$ at least $3$, i.e., for all $v \in T$, 
either $\deg(v) \geq 3$, or $\deg(v) = 1$.  
It is also assumed that the root $\rho$ is not a leaf, i.e., 
$\deg(\rho) \geq 3$.
\end{Assumption}


We give the following two equivalent formal definitions of ``topology''.
\begin{Definition} \label{def:topology}
\begin{itemize}
\item
Let $n$ be a positive integer and 
\begin{itemize}
\item $T$ be a tree with labeled leaves $v_1,\ldots,v_n$, so that $v_i$ is labeled by $i$.
\item $T'$ be a tree with labeled leaves $v'_1,\ldots,v'_n$, so that $v'_i$ is labeled by $i$.
\end{itemize}
We say that trees $T$ and $T'$ have the {\em same topology} if there exists a graph isomorphism $\phi : T \to T'$,
such that $\phi(v_i) = v'_i$, for $i=1,\ldots,n$.
\item {\em Equivalently}, the topology of $T$ is determined by the pairwise distances
$\left( d(v_i,v_i) \right)_{i,j=1}^n$, where $d$ is the graph-metric distance.
\item For a tree $T=(V,E)$ with labeled leaves $v_1,\ldots,v_n$, we write $\top(T)$ for the topology of $T$,
i.e. the equivalence class of $T$ under the relation of same topology (or $\top(T)$ is the
array of pairwise distances $\left( d(v_i,v_i) \right)_{i,j=1}^n$).
The topology of a tree $T$ rooted at $\rho$, is the topology
of the un-rooted tree $T$.
\end{itemize}
\end{Definition}

Naturally, it is impossible to reconstruct the topology with probability $1$.
\begin{Definition} \label{def:reconstruct}
Let $n$ be a positive integer and
\begin{itemize}
\item $\Trees$ be a family of rooted trees on $n$ labeled leaves.
\item $\Mut$ be a set of $|\A| \times |\A|$ stochastic matrices.
\end{itemize}
We write $\Trees \otimes \Mut$ for
\[
\Trees \otimes \Mut = \{ (T,(M^e)_{e \in E}) : T \in \Trees \mbox{ and } 
 \forall e \in E(T), M^e \in \Mut \}.
\]
We say that it is possible to reconstruct the topology from $k$ samples with probability $1 - \delta$,
if there exists a map 
$\psi : (|\A|^n)^k \to top(\Trees) = \{top(T) : T \in \Trees\}$,
such that for all $(T,(M^e)_e) \in \Trees \otimes \Mut$, if 
$(\sigma^t_{\partial T})_{t=1}^k$ are $k$ independent samples at the leaves, then
\[
\P\left[ \psi \left( (\sigma^t_{\partial T})_{t=1}^k \right) = \top(T) \right] \geq 1 - \delta.
\]
In this case we say that $\psi$ {\em reconstructs the topology} for $\Trees \otimes \Mut$ (from $k$
samples with probability $1 - \delta$). See (\ref{eq:diagram}) for a diagram representing $\psi$.
\begin{equation} \label{eq:diagram}
\begin{CD}
T \otimes M @>>> \sigma \\
@AA{\psi}A  @VV{\otimes^k}V \\
(\sigma^t_{\partial})_{t=1}^k @<<< (\sigma^t)_{t=1}^k
\end{CD}
\end{equation}
\end{Definition}
Note that this is a strong definition of reconstruction. In particular, if $\psi$ satisfies Definition
\ref{def:reconstruct}, then for any distribution of trees and matrices which is supported on $\Trees \otimes \Mut$,
$\psi$ reconstructs the underlying topology with probability at least $1 - \delta$. Also note that in
applications it is desirable that $\psi$ will have a simple algorithmic implementation.

\subsection{A Conjecture}
Our results are motivated by the following fundamental conjecture.

\begin{Conjecture} \label{conj:main}
Assume that the mutation matrices have a single order parameter
$\theta$, and let $\theta(e)$ be the order parameter for the mutation
matrix $M^e$ of edge $e$. 
Consider a markov random field on the $b+1$ regular tree where
the mutation matrices on all edges have the same parameter $\theta$.
Suppose that there exists $\theta_c$, such that 
\begin{itemize}
\item
If $\theta > \theta_c$, then the markov random field is in an ordered phase  
(in some technical sense, see below), 
\item
If $\theta < \theta_c$, then 
the markov random field is in an unordered phase.
\end{itemize}

We conjecture that the minimal number of samples needed in order to 
reconstruct phylogenies for the family of all trees on $n$ leaves, 
where all internal degree are at least $b+1$, is

\begin{itemize}
\item
$k = (c(\theta) + o(1)) \log n$, if for all edges of the phylogenetic tree, 
$\theta(e) \geq \theta > \theta_c$.
\item
$k = n^{c(\theta) + o(1)}$, if for all edges of the phylogenetic tree,
$\theta(e) \leq \theta < \theta_c$.
\end{itemize}
\end{Conjecture}
In the above conjecture the desired reconstruction probability is 
$1 - \delta$, for some fixed $0 < \delta < 1$. 
A more formal conjecture will require to specify how ``order'' is measured.
One possibility is to call a measure ordered for a specific value of $\theta$,
if the free measure on the infinite tree is extremal
(see subsection \ref{subsec:phase}). 
Another, is to look at spectral parameters of the
mutation matrices, and let $\theta(M) = |\lam_2(M)|$, where $|\lam_2(M)|$ is
the second largest eigen-value of $M$ in absolute value.
In this case it is natural to define $\theta_c$ by $b \theta_c^2 = 1$,
see \cite{KS,MP,KMP,JM}.
Following the results reported here, further support for 
Conjecture \ref{conj:main} was found in 
\cite{M:new,MS}

\subsection{The CFN model}
Below we focus on the model where $\Mut$ consists of all $2 \times 2$ matrices of the form
\begin{equation} \label{eq:CFN}
M^e = \left( \begin{array}{ll} 1 - \eps(e) & \eps(e) \\ \eps(e) & 1 - \eps(e) \end{array} \right),
\end{equation}
where $0 < \eps(e) < 1/2$ for all $e$.
We find it useful to denote $\theta(e) = 1 - 2 \eps(e)$.
Without loss of generality, we name the two colors $-1$ and $1$.
This model is referred to as the {\em Cavender-Farris-Neyman (CFN)} model.
It was studied \cite{Neyman,Farris,Cavender}, where it
is shown that if for all $e$ it holds that $\eps < \theta(e) < 1 - \eps$, then the
underlying topology can be reconstructed with probability $1 - \delta$ using $k = \poly_{\eps,\delta}(n)$ samples.
In \cite{St94} this result is generalized to mutation processes on 
any number of colors, 
provided that $M^e$ satisfy 
$\det(M^e) \notin [-1,-1+\eps] \cup [-\eps,\eps] \cup [1-\eps,1]$ 
for all $e$.
The dependency of $k$ on $\delta$ and $\eps$ is not stated explicitly
in these results.

It is desirable to minimize the number of samples needed for reconstruction.
Since the number of trees with $n$ leaves is exponential in 
$\Theta(n \log n)$, and each sample
consists of $n$ bits, 
it is clear that $\Omega(\log n)$ is a lower bound for the number
of samples.

In \cite{ESSW1,ESSW2} it is shown that for the CFN model (\ref{eq:CFN}) if for all $e$,
$1 > \theta_{\max} > \theta(e) > \theta_{\min} > 0$,
then it is possible to reconstruct the tree $T$ with probability
$1 - \delta$, if
\begin{equation} \label{eq:ESSW}
k > \frac{c \log n}{(1 - \theta_{\max})^2 \theta_{\min}^{d(T)}},
\end{equation}
where $d(T) = \Theta(\mbox{depth of } T)$ and $c = c(\delta)$.

For many of the trees that occur naturally in the reconstruction
setting, the depth of the tree is $\Theta(\log n)$. 
Bound (\ref{eq:ESSW}) on $k$ is therefore $k = n^{O(1)}$, which doesn't
improve previous bounds.
On the other hand, looking at families of random trees, $d(T)$ is typically $O(\log \log n)$, and therefore
by (\ref{eq:ESSW}) a $k = \polylog(n)$ number of samples suffice for reconstruction of a typical member of these families.

\subsection{Phase transition for the CFN model}
This paper is motivated by the following problem: 
When is the number of samples needed in order to reconstruct the
topology of $T$ polynomial in $n$ and when is it poly-logarithmic
in $n$? 
The hardest case in the analysis of \cite{ESSW1,ESSW2}
is that of a balanced tree. We will focus on balanced trees below.

\begin{Definition} \label{def:balanced}
A tree $T$ rooted at $\rho$ is {\em balanced}, if all the leaves of
$T$ have the same distance to $\rho$, i.e., there exists an $r$ such that
\[
\partial T = \{v \in V(T) : d(v,\rho) = r\}.
\] 
\end{Definition}

We first focus on the case where the mutation rate is 
the same for all edges.
The following two theorems already indicate
the importance of certain ``phase-transitions'' for the problem.
For $b \geq 2$, we let $\Trees_b^{\ast}((b+1) b^q)$ 
denote the space of all balanced rooted trees
on $n = (b+1) b^q$ leaves, where all the internal degrees are exactly
$b+1$. We call a tree in $\Trees_b^{\ast}((b+1) b^q)$, a 
{\em $(q+1)$-level $(b+1)$-{\em regular} tree}.

\begin{Theorem} \label{thm:CFN1}

Consider the tree reconstruction problem for the CFN model
on the space $\Trees_b^{\ast}(n)$, where $n = (b+1) b^q$, and  
for all $e$, $\theta(e) = \theta$, is independent of $e$.
If $b \theta^2 > 1$, then there exists $c_{\theta} < \infty$ such that
for all $\delta > 0$, it is possible to
reconstruct the topology from $k$ samples with probability $1 -
\delta$, where $k = c_{\theta} (\log n - \log \delta)$.
\end{Theorem}

This result could not be extended to $\theta$ such that $b \theta^2 < 1$,
as the following theorem implies that reconstructing
balanced trees actually requires polynomial number of samples 
when the mutation rate is high.
\begin{Theorem} \label{thm:CFN_lower}
Suppose that $b \theta^2 < 1$, then there exists $q_0$ such that for
all $q \geq q_0$, the tree reconstruction problem for the CFN model
on the space $\Trees_b^{\ast}(n)$, where $n = (b+1) b^q$, $q \geq q_0$
and for all $e$, $\theta(e) = \theta$, satisfies the following.
 
Given a uniformly chosen tree from $\Trees_b^{\ast}((b+1) b^q)$
(assume that the initial distribution of the
color at the root is uniform, $\pm 1$ with probability $1/2$ each)
and $k$ samples of the coloring at the boundary of the tree,  
the probability of reconstructing the topology is at most
\begin{equation} \label{eq:CFN_lower}
k (b \theta^2)^{q - \log_b q - \log_b (-\log b \theta^2)} =
O \left(k n^{(1 + 2 \log_b \theta)(1 - \frac{\log_b \log_b n}{\log_b
      n}) } \right).
\end{equation}
\end{Theorem}

Theorems \ref{thm:CFN1} and \ref{thm:CFN_lower} indicate the 
importance of the study of the phase
transitions for the Ising model on trees, 
where an interesting phase transition occurs when
$b \theta^2 = 1$, see Subsection \ref{subsec:phase} for more background.

Later we generalize Theorem \ref{thm:CFN1} to the standard model 
on balanced trees. Let $\Trees_{\geq b}^{\ast}(n)$ be the space of all
balanced rooted trees on $n$ leaves,
where all the internal degrees are at least $b+1$.

\begin{Theorem} \label{thm:CFN2}
Consider the tree reconstruction problem for the CFN model on the
space $\Trees_{\geq b}^{\ast}(n)$.
Suppose that
$\theta_{\min}$ satisfies $b \theta_{\min}^2 > 1$, and that all edges $e$ satisfy that 
$\theta_{\min} \leq \theta(e) \leq \theta_{\max} < 1$.
Then there exists a constant 
$c = c(\theta_{\min},\theta_{\max}) = c(\theta_{\min})/(1 -
\theta_{\max})^2 < \infty$, such that for all $\delta > 0$, it is possible to
reconstruct the topology with probability $1 - \delta$ from 
$k = c (\log n - \log \delta) $ samples in
$\poly_{\delta,\theta_{\min},\theta_{\max}}(n)$ time.
\end{Theorem}

Theorem \ref{thm:CFN2} implies in particular a conjecture of Steel \cite{Stconj} for balanced trees,
which initiated this work. We believe the Theorem \ref{thm:CFN2} could play an important role in proving the analogous
result for general (non-balanced) trees.
We can also prove an upper bound for the number of samples when 
$b \theta_{\min}^2 < 1$.

\begin{Theorem} \label{thm:CFN3}
Consider the tree reconstruction problem for the CFN model on the
space $\Trees_{\geq b}^{\ast}(n)$, and let $q = \log_b n$.
Suppose $\theta_{\max} < 1$, 
$\theta_{\min}$ satisfies $g^2 < b \theta_{\min}^2 < 1$, and all edges
$e$ satisfy $\theta_{\min} \leq \theta(e) \leq \theta_{\max}$.
Then for all $\delta > 0$, 
it is possible to reconstruct the topology with probability
$1 - \delta$ given 
$k = c(\theta_{\min},\theta_{\max}, \delta,g)\,\, g^{- 8q} $ samples in
$\poly_{\delta,\theta_{\min},\theta_{\max}}(n)$ time.
\end{Theorem}

Theorem \ref{thm:CFN_lower} also implies a lower bound on learning the
tree in the PAC setting, see \cite{ADFK,CGG}.


\subsection{Phase transitions for the Ising model on the tree} \label{subsec:phase}
The extremality of the free measure for the Ising model on the regular tree plays a crucial
role in this paper. 
The study of the extremality of the free measure begins
with \cite{Sp,Hi}. 
Latter papers include \cite{BRZ,EKPS,Io,M:recur}
(see \cite{EKPS} for more detailed background).

Consider the CFN model on a $q$-level $(b+1)$-regular tree, 
where $\theta(e) = \theta$ for all $e$,
and where the root is chosen to be each of the two colors 
with probability $1/2$.
This measure is known in statistical physics as the free Gibbs
 measure for the Ising model on the homogeneous tree (or
Bethe lattice). Note in particular, that the tree topology is fixed in advance.
Given a {\em single} sample of the colors at the leaves of the tree, we want to reconstruct
some information on the root color. The basic question is whether 
amount of information that can be
reconstructed decay to $0$, as $q$ increases. Let $\sigma_{\rho}$
denote the color of the root, and $\sigma_q$ denote
the colors at level $q$. It turns out that the following conditions are equivalent.
\begin{itemize}
\item
$I(\sigma_{\rho},\sigma_q) \to 0$, where $I$ is the mutual information operator.
\item
The total variation distance between the distribution of $\sigma_q$
given $\sigma_{\rho}=1$, and the distribution of $\sigma_q$ given $\sigma_{\rho}=-1$
decays to $0$ as $q \to \infty$.
\item
For all algorithms, the probability of reconstructing $\sigma_{\rho}$ from $\sigma_q$ decays to $1/2$
as $q \to \infty$.
\item
The free Gibbs measure for the Ising model on the
infinite $(b+1)$-regular tree is extremal.
\end{itemize}
(see \cite{EKPS} for definitions and proof of the equivalence,
and \cite{M:lam2} for this equivalence for general Markov random 
fields on the tree).

The extremality phase transition may be formulated as follows. 
When $b \theta^2 > 1$, 
some information on the root can be reconstructed independently
of the height of the tree, i.e., none of the equivalent conditions above hold.
When $b \theta^2 < 1$, 
all of the above conditions hold, and it is therefore impossible to reconstruct
the root color, as $q \to \infty$.

Theorem \ref{thm:CFN1} is based on algorithmic aspects of this phase transition
discussed in \cite{M:recur}, 
while Theorem \ref{thm:CFN_lower} utilizes information bounds from \cite{EKPS}.




\subsection{Paper outline}
In Section \ref{sec:fixed} we give a short proof of Theorem \ref{thm:CFN1} 
as it demonstrates some of the key ideas to be applied later in
Theorem \ref{thm:CFN2} and Theorem \ref{thm:CFN3}.
In Section \ref{sec:lower_bound} we prove Theorem
\ref{thm:CFN_lower}. The proof uses information
bounds and is somewhat independent from the other sections.
In Section \ref{sec:rec_maj} we study in detail the behavior of
majority algorithms for local reconstruction
as the main technical ingredient to be used later.
Section \ref{sec:four_points} contains some basic results regarding
large deviations and four-point conditions.
In section \ref{sec:alg} we present 
the proofs of theorems \ref{thm:CFN2} and \ref{thm:CFN3}.

\section{Logarithmic reconstruction for fixed $\theta$} \label{sec:fixed}
We start by proving Theorem \ref{thm:CFN1}.
We first define formally the function $\Maj$. Note that when the number of inputs
is even, this function is {\em randomized}.
\begin{Definition} \label{def:maj}
Let $\Maj : \{-1,1\}^d \to \{-1,1\}$ be defined as:
\[
\Maj(x_1,\ldots,x_d) = \sign(\sum_{i=1}^d x_i + 0.5 \omega),
\]
where $\omega$ is an unbiased $\pm 1$ variable which is independent of the $x_i$.
Thus when $d$ is odd,
\[\Maj(x_1,\ldots,x_d) = \sign(\sum_{i=1}^d x_i).
\]
When $d$ is even,
\[
\Maj(x_1,\ldots,x_d) = \sign(\sum_{i=1}^d x_i),
\]
unless $\sum_{i=1}^d x_i = 0$, in which case $\Maj(x_1,\ldots,x_d)$ is chosen to be $\pm 1$ with probability $1/2$.
\end{Definition}

For $b \geq 2$, we call a tree $T$ rooted at $\rho$, 
the {\em $\ell$ level $b$-ary} tree, if
all internal nodes have exactly $b$ descendants and all the leaves are
at distance $\ell$ from the root.

\begin{Lemma} \label{lem:recur_fix_eps}
Let $b$ and $\theta$ be such that $b \theta^2 > 1$.
Then there exists $\ell = \ell(\theta)$, 
and $1 > \eta_0 = \eta_0(\ell,\theta) > 0$,
such that for all $\eta \geq \eta_0$, 
the CFN model on the $\ell$-level $b$-ary tree $T$ with
\begin{itemize}
\item $\theta(e) = \theta$ for all $e$ which is not adjacent to
  $\partial T$, and
\item $\theta(e) = \theta \eta$ for all $e$ which is adjacent to $\partial T$.
\end{itemize}
satisfies
\[
\E[+\Maj(\sigma_{\partial T})| \sigma_{\rho} = +1] =
\E[-\Maj(\sigma_{\partial T})| \sigma_{\rho} = -1] \geq \eta_0.
\]
\end{Lemma}
This follows from Section 3 of \cite{M:recur}.
Lemma \ref{lem:recur_fix_eps} also follows from the more general  
Theorem \ref{thm:maj_good} below.

The only other tool needed for the proof in this case are standard 
large deviations results, see, e.g., \cite[Corollary A.1.7]{AS}.
\begin{Lemma} \label{lem:dev_fix_eps}
Let $S = \sum_{i=1}^k X_i$, where $X_i$ are i.i.d. $\{-1,1\}$ random variables.
Then for all $a > 0$,
\[
\P \left[ |S - \E[S]| \geq a \right] \leq 2 \exp(- \frac{a^2}{2 k}).
\]
\end{Lemma}

\begin{Definition} \label{def:ast}
Let $T$ be a balanced tree.
\begin{itemize}
\item
The {\em $\ell$-topology of $T$}
is the function 
$d^{\ast}_{\ell} : \partial T \times \partial T \to \{0,\ldots, 2 \ell + 2
\}$, defined by $d^{\ast}_{\ell}(u,v) = \min \{d(u,v), 2 \ell + 2\}$.
\item 
We let $L_{\partial - i} = \{ v \in T : d(v,\partial T) = i\}$.
\item
The $\ell$ labeling of $T$ is the
labeling of $\cup_{i=0}^{\ell} L_{\partial - i}$, where 
$v \in L_{\partial - i}$ is
labeled by 
\[
\partial T(v) = \{w \in \partial T : d(v,w) = i\}.
\]
\end{itemize}
\end{Definition}
Note that for a balanced tree $T$, the $\ell$-topology of $T$
determines the $\ell$-labeling of $T$ -- for $i \leq \ell$
the labels of $L_{\partial - i}$ are given by the sets 
$\left\{ \{ w' \in \partial T: d^{\ast}_{\ell}(w,w') \leq 2 i\} : 
w \in \partial T
\right\}$. Moreover, if $u,v \in V$, $d(u,\partial T) \leq \ell$ 
and $d(v, \partial T) \leq \ell$, then $v$ is a descendant of $u$ iff
$\partial T(v) \subset \partial T(u)$.

The core of the proof of Theorem \ref{thm:CFN1} is the following lemma.
\begin{Lemma} \label{lem:rec_fix_eps}
Let $b$ and $\theta$ be such that $b \theta^2 > 1$.
Let $\ell$ and $\eta_0$ be such that Lemma \ref{lem:recur_fix_eps} holds,
and assume that $\eta \geq \eta_0$.
Consider the CFN model on the family of balanced tree of $q$ levels, where all
internal nodes have at least $b$ children and the total number of
leaves is $n$. Assume that $\theta : E \to [0,1]$ satisfies 
\begin{itemize}
\item $\theta(e) = \theta$ for all $e$ which is not adjacent to 
$\partial T$, and
\item $\theta(e) = \theta \eta$ for all $e$ which is adjacent to $\partial T$.
\end{itemize}
Then
\begin{itemize}
\item
given $k$ independent samples of the process at the leaves of $T$,
$(\sigma^t_{\partial T})_{t=1}^k$, and $\ell \leq q$,
it is possible to recover the $\ell$-topology of $T$ 
with error probability bounded by
\begin{equation} \label{eq:bnd_fix}
n^2 \exp(-c^{\ast}\, k),
\end{equation}
where $c^{\ast} = \eta_0^4 \theta^{4 \ell} (1 - \theta^2)^2/8$.
\item
For all $T$ and $i \geq 0$, there exists a map 
$\Psi = \Psi_T  
: \{\pm 1\}^{\partial T} \to \{\pm 1\}^{L_{\partial - i \ell}}$ for
which the following hold.

If $\sigma$ is distributed according to the CFN
model on $T$, and $\sigma' = \Psi(\sigma_{\partial T})$, then
$(\sigma'_v)_{v \in L_{\partial - i \ell}} = 
(\sigma_v \tau_v)_{v \in L_{\partial - i \ell}}$, where
$\tau_v$ are i.i.d. variables.
Moreover, $\tau_v$ are  
independent of $(\sigma_v)_{d(v,\partial T) \geq i \ell}$ and satisfy
$\E[\tau_v] \geq \eta_0$.

The map $\Psi$ may be constructed from the $(i \ell)$-topology of $T$.
In particular, if $T_1$ and $T_2$ have the same $(i \ell)$-topology, then 
$\Psi_{T_1} = \Psi_{T_2}$.

\end{itemize}
\end{Lemma}

\proofs

Let $c(u,v)$ be the correlation between $u$ and $v$
\[
c(u,v) = \frac{1}{k}\sum_{t=1}^k \sigma^t_u \sigma^t_v.
\]
Suppose that $d(u,v) = 2r$.
Then $\E[c(u,v)] = \alpha_r$, where $\alpha_r = \eta^2 \theta^{2 r}$.
We let
\[
I_r = \left\{ \begin{array}{ll}
          \left(\frac{\alpha_{r+1} + \alpha_r}{2},\frac{\alpha_r + \alpha_{r-1}}{2} \right) &
          \mbox{ if } 1 \leq r \leq \ell, \\
          \left[-1,\frac{\alpha_{\ell+1} + \alpha_{\ell}}{2} \right) &
          \mbox{ if } r = \ell+1. \end{array} \right.
\]
Since $k c(u,v)$ is a sum of $k$ i.i.d. $\pm 1$ variables, it follows from Lemma \ref{lem:dev_fix_eps} that
for all $u$ and $v$,
\begin{eqnarray*}
\P[c(u,v) \notin I_{d(u,v)/2}] &\leq&
\max_{1 \leq r \leq \ell} 2 \exp 
\left(-\frac{1}{2 k} \left(k \frac{\alpha_{r+1} - \alpha_{r}}{2}
\right)^2 \right) \\ &=&
2 \exp \left(-k \eta^4 \theta^{4 \ell} (1 - \theta^2)^2/8 \right) \leq
2 \exp \left(-k \eta_0^4 \theta^{4 \ell} (1 - \theta^2)^2/8 \right).
\end{eqnarray*}
Note that the intervals $(I_r)_{r=1}^{\ell+1}$ are disjoint. 
Define $D^{\ast}_{\ell}(u,v) = 2 r$,
if $c(u,v) \in I_r$.
Then $D^{\ast}_{\ell}(u,v) = d^{\ast}_{\ell}(u,v)$, for all $u$ and $v$,
with error probability bounded by (\ref{eq:bnd_fix}), 
thus proving the first claim of the lemma.

We now prove the second claim by induction on $i$. The claim is
trivial for $i = 0$ as we may take $\Psi$ to be the identity map.

For the induction step, suppose that we are given $d^{\ast}_{i \ell+\ell}$.
Label all the vertices $v \in \cup_{j=0}^{(i+1) \ell} L_{\partial - j}$ by
the $(i+1) \ell$-labeling of $T$. 
By the induction hypothesis, there exists a map $\Psi'$ such that 
$\Psi'(\sigma_{\partial T}) = 
(\sigma_v \tau_v)_{v \in L_{\partial T - i \ell}}$, 
where $\tau_v$ are i.i.d. ${\pm 1}$ variables.
Moreover, $\tau_v$ are  independent of 
$(\sigma_v)_{d(v,\partial T) \geq i \ell}$ and satisfy
 $\E[\tau_v] \geq \eta_0$.  

By the properties of the labeling, for each 
$w \in L_{\partial - (i+1)\ell}$, there exists a set 
$R(w) \subset L_{\partial - i \ell}$, 
which is the set of leaves of an
$\ell$-level $b$-ary tree rooted at $w$.

We now let  $\Psi(\sigma_{\partial T}) = 
(\widehat{\sigma}_w)_{w \in L_{\partial - i \ell - \ell}}$, where
\[
\widehat{\sigma}_w = 
\Maj((\Psi'(\sigma_{\partial T}))_v : v \in R(w)) = 
\Maj(\sigma_v \tau_v : v \in R(w)).
\]
By Lemma \ref{lem:recur_fix_eps} it follows that 
$(\widehat{\sigma}_w)_{w \in L_{\partial - i \ell - \ell}} =  
(\sigma_w \tau_w)_{w \in L_{\partial - i \ell - \ell}}$,
where $\tau_w$ are i.i.d. $\pm 1$ variables. Moreover, 
$\tau_w$ are independent 
of $(\sigma_v)_{d(v,\partial T) \geq i \ell + \ell}$ and satisfy
$\E[\tau_w] \geq \eta_0$, proving the second claim.
$\QED$

\prooft{of Theorem \ref{thm:CFN1}}
Let $b$ and $\theta$ be such that $b \theta^2 > 1$.
Let $\ell$ and $\eta_0$ be such that Lemma \ref{lem:recur_fix_eps}
holds. Note that if $i \ell \geq q$, then $d^{\ast}_{i \ell} = d$.
Therefore, in order to recover $d$, it suffices to apply Lemma 
\ref{lem:rec_fix_eps} recursively in order to recover $d^{\ast}_{i \ell}$, 
for $i = 0, \ldots, \lceil q/\ell \rceil$. 

It is trivial to recover $d^{\ast}_0(v,u) = 2 \one_{v \neq u}$.
We now show how given $d^{\ast}_{i \ell}$ and the samples 
$(\sigma^t_\partial)_{t=1}^k$, we can recover $d^{\ast}_{i \ell+\ell}$ with
error probability bounded by $n^2 \exp (-c^{\ast} k) / b^{2 \ell}$. 

Let $\Psi_i : \{\pm 1\}^{\partial T} \to \{\pm 1\}^{L_{\partial - i
    \ell}}$ 
be the function defined in  
second part of Lemma \ref{lem:rec_fix_eps} given $d^{\ast}_{i \ell}$.

Then
\[ 
(\Psi_i(\sigma^t_{\partial T}))_{t=1}^k = 
(\sigma^t_v \tau^t_v : v \in L_{\partial - i \ell})_{t=1}^k,
\]
where $\tau^t_v$ are i.i.d. variables with $\E[\tau^t_v] \geq
\eta_0$. Moreover, $\tau^t_v$ are independent of 
$(\sigma^t_v : d(v,\partial T) \geq i \ell, 1 \leq t \leq k)$.

By the first part of the lemma, given
$(\sigma^t_v \tau^t_v : v \in L_{\partial - i \ell})_{t=1}^k,$
we may recover 
\[
d' : L_{\partial - i \ell} \times L_{\partial - i \ell} \to 
\{0,\ldots, 2 \ell + 2\},
\]
defined by $d'(u,v) = \min \{d(u,v), 2 \ell + 2\}$, 
with error probability bounded
by $n^2 \exp (-c^{\ast} k) / b^{2 i \ell}$.

Note that
\begin{equation} \label{eq:recd}
d^{\ast}_{i \ell + \ell}(u,v) = \left\{ \begin{array}{ll}
     d^{\ast}_{i \ell}(u,v) & \mbox{if } 
     d^{\ast}_{i \ell}(u,v) \leq 2 i \ell, \\
     d'(u',v') + 2 i \ell  & 
     \mbox{if } u \in \partial T(u'), v \in \partial T(v'), 
     \{u',v'\} \subset L_{\partial - i \ell}, u' \neq v'.
     \end{array} \right.
\end{equation}        
Thus given $d^{\ast}_{i \ell}$, by recovering $d'$, 
we may recover $d^{\ast}_{i \ell + \ell}$.

Let $A_i$ be the event of error in recovering
$d^{\ast}_{i \ell + \ell}$ given $d^{\ast}_{i \ell}$ and
$\alpha = \sum_{i=0}^{\lceil q/l \rceil} \P[A_i]$.
Then the probability of error in the recursive scheme above is bounded by 
$\alpha$ and
\[
\alpha \leq \exp (-c^{\ast} k)
\left(n^2 + n^2/b^{2 \ell} + n^2/b^{4 \ell} + \cdots \right) \leq
2 n^2 \exp (-c^{\ast} k).
\] 
Defining ${c'}^{-1} = c^{\ast}$, and taking
\[
k = \frac{\log(2 n^2) - \log \delta}{c^{\ast}} = 
c'(2 \log n + \log 2 - \log \delta)
\]
we obtain $\alpha \leq \delta$.
The statement of the theorem follows by letting $c_{\theta} = 3 c'$.
$\QED$

\section{Polynomial lower bound for $b \theta^2 < 1$} \label{sec:lower_bound}

In this section we prove Theorem \ref{thm:CFN_lower} via an entropy argument.
Let $X$ and $Y$ be discrete random variables. 
Recall the definitions of the entropy of $X$, $H(X)$, 
the conditional entropy of $X$
given $Y$, $H(X | Y)$, and the mutual information of $X$ and $Y$,
$I(X,Y)$:
\begin{eqnarray*}
H(X) &=& -\sum_x \P[X = x] \log_2 \P[X = x] , \\
H(X|Y) &=& \E_y H(X | Y = y) = H(X,Y) - H(Y), \\
I(X,Y) &=& H(X) + H(Y) - H(X,Y) = H(X) - H(X|Y) = H(Y) - H(Y|X).
\end{eqnarray*}
(see e.g. \cite{CT} for basic properties of $H$ and $I$).

The core of the proof of Theorem \ref{thm:CFN_lower} 
is the fact that for the $q$-level $b$-ary tree, if $b \theta^2 < 1$, 
then the correlation between the color at the root and the coloring
of the boundary of the tree decays exponentially in $q$.
We will utilize the following formulation from \cite{EKPS} 
(see also \cite{BRZ,Io}).

\begin{Lemma}[\cite{EKPS}] \label{lem:ekps_inf}
Let $\sigma$ be a sample of the CFN process on the $q$-level $b$-ary
tree. Then
\[
I(\sigma_{\rho},\sigma_{\partial T}) \leq b^q \theta^{2q}.
\]
\end{Lemma}

We will also use some basic properties of $I$, see e.g. \cite{CT}.
\begin{Lemma} \label{lem:dpl}
Let $X,Y$ and $Z$ be random variables such that $X$ and $Z$ are independent given $Y$, then
\begin{equation} \label{eq:dpl}
I(X,Z) \leq \min\{I(X,Y),I(Y,Z)\} \,\,\,\,\,
\mbox{(``Data Processing Lemma'')},
\end{equation}
\begin{equation} \label{eq:Itr}
I((X,Y),Z) = I(Y,Z),
\end{equation}
\begin{equation} \label{eq:cond_ind}
I((X,Z),Y) \leq I(X,Y) + I(Z,Y).
\end{equation}
\end{Lemma}

It is well known that if $I(X,Y)$ is small, 
then it is hard to reconstruct $X$ given $Y$.

\begin{Lemma}[Fano's inequality] \label{lem:fano}
Let $X$ and $Y$ be random variables s.t.
$X$ takes values in a set $A$ of size $m$, $Y$ takes values
in a set $B$, and
\begin{equation} \label{eq:defDelta}
\Delta = \Delta(X,Y) = \sup_{f : B \to A} \P[f(Y) = X],
\end{equation}
is the probability of reconstructing the value of $X$ given $Y$
(the $\sup$ is taken over all randomized functions).
Then
\begin{equation} \label{eq:inf_and_rec}
H(\Delta) + (1 - \Delta) \log_2 (m - 1) \geq H(X|Y),
\end{equation}
where $H(\Delta) = - \Delta \log_2 \Delta - (1 - \Delta) \log_2 (1 - \Delta)$.
\end{Lemma}

It is helpful to have the following easy formula.
\begin{Lemma} \label{lem:ntop}
The number of topologies for $\ell$ level $b$-ary trees 
on $b^{\ell}$ labeled leaves, $n_{top}(\ell)$, is
\begin{equation} \label{eq:n_topologies}
n_{top}(\ell) = \frac{b^{\ell}!}{b!^{\sum_{j=0}^{\ell-1} b^j}}.
\end{equation}
In particular, if $\ell \geq b^3$, then
\begin{equation} \label{eq:log_n_topologies}
\log n_{top}(\ell) \geq b^{\ell-1} \log (b^{\ell}).
\end{equation}
\end{Lemma}

\proofs
Clearly $n_{top}(1) = 1$, and
\[
n_{top}(\ell) = \frac{n_{top}(\ell-1)^b}{b!} \binom{b^{\ell}}{b^{\ell-1}\,\cdots\,b^{\ell-1}}.
\]
We therefore obtain (\ref{eq:n_topologies}) by induction.
To obtain (\ref{eq:log_n_topologies}), we note that by Stirling's formula,
for $\ell \geq b^3$,
\begin{eqnarray*}
\log \left( n_{top}(\ell) \right) &=& \log (b^{\ell}!) - 
\log (b!) \sum_{j=0}^{\ell-1} b^j \geq
b^{\ell} \log (b^{\ell}) - b^{\ell} - 
\frac{b^{\ell} - 1}{b - 1} \log(b!) \\ &\geq&
b^{\ell} \log(b^{\ell}) \left( 1 - \frac{1}{\log(b^{\ell})} - 
               \frac{1}{b^2} \right)  \geq
b^{\ell - 1} \log (b^{\ell}).
\end{eqnarray*}
$\QED$

\prooft{of Theorem \ref{thm:CFN_lower}}
We first note that given the topology 
of a tree $T \in \Trees^{\ast}_b((b+1) b^q)$, the root of $T$ is
uniquely determined as the unique vertex which has the same distance
to all the leaves. Therefore if 
$T,T' \in \Trees^{\ast}_b((b+1) b^q)$ have the same
topology, then $T$ and $T'$ are isomorphic as rooted trees, i.e.,
there exists a graph homomorphism $\psi$ of $T$ onto $T'$ which maps
leaves to leaves of the same label and the root of $T$ to the
root of $T'$. In the proof below we won't distinguish between the
topology of $T$ and $T$.
 
Assuming $b \theta^2 < 1$, 
we want to prove a lower bound on the number of samples needed in order
to reconstruct the tree, given that the tree is chosen uniformly at random.
Clearly, the probability of reconstruction increases, if in addition to the
samples we are given additional information.
We will assume that we are given the $(q - \ell + 1)$-topology of tree, 
i.e., for all $u,v \in \partial T$, we are given
$d^{\ast}(u,v) = \min\{ d(u,v), 2 (q - \ell) + 4\}$; the value of $\ell < q$
will be specified later.

Given $d^{\ast}$, we have the $(q-\ell+1)$-labeling of
the nodes of $L_{\ell} =  \{ v: d(v,\rho) = \ell\}$.
Therefore, the reconstruction problem
reduces to reconstructing an element of $T_b^{\ast}((b+1) b^{\ell-1})$
on the set of labeled leaves $\{ \partial T(v) : v \in L_{\ell}\}$.
Moreover, it is easy to see that given $d^{\ast}$, the conditional
distribution over $T_b^{\ast}((b+1) b^{\ell-1})$ is uniform.

Recall that $\sigma^t_{\partial} = (\sigma^t_v : v \in \partial T)$.
We may generate $(\sigma^t_{\partial})_{t=1}^k$, 
in the following manner:
\begin{itemize}
\item Choose $T \in \Trees_b^{\ast}((b+1) b^{\ell-1})$ uniformly at random.
\item Given $T$, generate the samples at level
$l$, i.e, $\sigma^t_{\ell} = (\sigma^t_v : v \in L_{\ell})$, 
for $1 \leq t \leq k$.
\item Given $(\sigma^t_{\ell} : 1 \leq t \leq k)$, generate 
$(\sigma^t_{\partial})_{t=1}^k$.
\end{itemize}

We conclude
that $T$ and $(\sigma^t_{\partial})_{t=1}^k$ are conditionally independent
given $(\sigma^t_{\ell})_{t=1}^k$.

In particular, by the Data processing Lemma (\ref{eq:dpl}),
\begin{equation} \label{eq:dpl_top}
I \left(T,(\sigma^t_{\partial})_{t=1}^k \right) \leq
I \left( (\sigma^t_{\ell})_{t=1}^k,(\sigma^t_{\partial})_{t=1}^k \right).
\end{equation}

Since $\sigma^t$ and $\sigma^{t'}$ are independent for $t \neq t'$,
\begin{equation} \label{eq:inf1}
I \left( (\sigma^t_{\ell})_{t=1}^k,(\sigma^t_{\partial})_{t=1}^k \right) =
\sum_{t=1}^k I(\sigma^t_{\ell},\sigma^t_{\partial}).
\end{equation}

Let $\sigma^t_{\partial T(v)}$ be the coloring of the set $\partial T(v)$.
Note that $(\sigma^t_{\partial T(v)} : v \in L_{\ell})$ are conditionally independent given $\sigma^t_{\ell}$.
Therefore, by (\ref{eq:cond_ind}),
\begin{equation} \label{eq:inf2}
I(\sigma^t_{\ell},\sigma^t_{\partial}) \leq \sum_{v \in L_{\ell}}
I(\sigma^t_{\ell},\sigma^t_{\partial T(v)}).
\end{equation}
Finally note that $(\sigma^t_w : w \in L_{\ell}, w \neq v)$ are
independent of $\sigma^t_{\partial T(v)}$
given $\sigma_v^t$, and therefore for all $v \in L_{\ell}$,
\begin{equation} \label{eq:inf3}
I(\sigma^t_{\ell},\sigma^t_{\partial T(v)}) =
I(\sigma_v^t,\sigma^t_{\partial T(v)}).
\end{equation}

Combining (\ref{eq:dpl_top}), (\ref{eq:inf1}), (\ref{eq:inf2}) and (\ref{eq:inf3}) we obtain
\[
I \left(T,(\sigma^t_{\partial})_{t=1}^k \right) \leq
\sum_{t=1}^k \sum_{v \in L_{\ell}} I(\sigma_v^t,\sigma^t_{\partial T(v)}).
\]
By Lemma \ref{lem:ekps_inf}, $I(\sigma_v^t,\sigma^t_{\partial T(v)})
\leq b^{q - \ell + 1} \theta^{2(q - \ell + 1)} \leq 
b^{q - \ell} \theta^{2(q - \ell)}$.
Therefore
\[
I \left(T,(\sigma^t_{\partial})_{t=1}^k \right) \leq 
k (b+1) b^q \theta^{2(q - \ell)}.
\]
Letting $m'$ be the number of topologies for trees in 
$\Trees^{\ast}((b+1) b^{\ell-1})$, we see that
\[
H \left(T \,|\, (\sigma^t_{\partial})_{t=1}^k \right) = H(T) - I \left(T,(\sigma^t_{\partial})_{t=1}^k \right) \geq
\log_2 m' - k (b+1) b^q \theta^{2(q - \ell)}.
\]

By lemma \ref{lem:fano}, we conclude that the probability
$\Delta = \Delta \left( T, (\sigma^t_{\partial})_{t=1}^k \right)$, of reconstructing
$T$ given $(\sigma^t_{\partial})_{t=1}^k$, satisfies
\begin{equation} \label{eq:fano_app}
H(\Delta) + (1 - \Delta) \log_2 m' \geq \log_2 m' - 
k (b + 1) b^q \theta^{2(q - \ell)},
\end{equation}
The rest of the proof consists of calculations showing
how to derive (\ref{eq:CFN_lower}) from (\ref{eq:fano_app}).

Clearly $m'$ is at least $m = n_{top}(\ell)$.
Rewriting (\ref{eq:fano_app}), we obtain
\[
H(\Delta) + k (b + 1) b^q \theta^{2(q - \ell)} \geq \Delta \log_2 m' \geq 
\Delta \log_2 m,
\]
from which we conclude that
\begin{equation} \label{eq:fin2_lower}
\Delta \leq \max \left\{\frac{2 H(\Delta)}{\log_2 m}, 
\frac{2 k (b + 1) b^q \theta^{2(q - \ell)}}{\log_2 m} \right\}.
\end{equation}

Note that $- (1-x) \log (1-x) \leq x$ for $x \in [0,1]$, and therefore
$H(\Delta) \leq - \Delta \log_2 \Delta + \Delta \log_2(e)$.
Thus if
$
\Delta \leq {2 H(\Delta)}/{\log_2 m},
$
then
$
0.5 \Delta \log_2(m) \leq -\Delta \log_2 \Delta + \Delta \log_2(e),
$
or $\Delta \leq e / \sqrt{m}$.
So by (\ref{eq:fin2_lower}), we obtain
\begin{equation} \label{eq:fin3_lower}
\Delta \leq \max \left\{ \frac{e}{\sqrt{m}},
\frac{2 k (b + 1) b^q \theta^{2(q - \ell)}}{\log_2 m} \right\}.
\end{equation}
Therefore, if $\ell \geq b^5$ (say) then by Lemma \ref{lem:ntop},
\[
\Delta \leq \max \left\{
\exp \left(1 - 0.5 b^{\ell - 1} \log b^{\ell} \right),
\frac{2 k b (b + 1)(b \theta^2)^{q - \ell}}{\log_2 b^{\ell}} \right\}
\leq
\max\{ \exp(-b^{\ell+1}), k (b \theta^2)^{q - \ell} \}.
\]
We now take $\ell = \lfloor \log_b q + \log_b (-\log b \theta^2) \rfloor$, so
$\exp(-b^{\ell+1}) \leq (b \theta^2)^q$.
Since we have the freedom of choosing $\ell$, we conclude that
\[
\Delta \leq k (b \theta^2)^{q - \log_b q - \log_b (-\log b \theta^2)},
\]
for large $q$ as needed.
$\QED$

\section{Majority on trees} \label{sec:rec_maj}

In this section
we analyze the behavior of the majority algorithm on balanced $b$-ary
trees. Theorem \ref{thm:maj_good} will be used later in the proof of
theorems \ref{thm:CFN2} and \ref{thm:CFN3}

\begin{Definition} \label{def:CFN_t_e}
Let $T=(V,E)$ be a tree rooted at $\rho$ with boundary 
$\partial T$.
For functions $\theta' : E \to [0,1]$ and $\eta' : \partial T \to [0,1]$, let $CFN(\theta',\eta')$
be the CFN model on $T$ where
\begin{itemize}
\item $\theta(e) = \theta'(e)$ for all $e$ which is not adjacent to $\partial T$, and
\item $\theta(e) = \theta'(e) \eta'(v)$ for all $e=(u,v)$, 
with $v \in \partial T$.
\end{itemize}
Let
\[
\widehat{\Maj}(\theta',\eta') = \E[+\Maj(\sigma_{\partial T})| \sigma_{\rho} = +1] =
\E[-\Maj(\sigma_{\partial T})| \sigma_{\rho} = -1],
\]
where $\sigma$ is drawn according to $CFN(\theta',\eta')$.
\end{Definition}
For functions $\theta$ and $\eta$ as above we'll abbreviate by writing $\min \theta$ for $\min_{E} \theta(e)$,
$\max \eta$ for $\max_{v \in \partial T} \eta(v)$, etc.
The function $\widehat{\Maj}$ measures how well the majority
calculates the color at the root of the tree.

\begin{Theorem} \label{thm:maj_good}
Let 
\begin{equation} \label{eq:def_a}
a(d) = 2^{1-d} \lceil \frac{d}{2} \rceil \binom{d}{\lceil \frac{d}{2} \rceil}.
\end{equation}
For all $\ell$ integer, $\theta_{\min} \in [0,1]$ and 
$0 \leq \alpha < a(b^{\ell}) \theta_{\min}^{\ell}$, there exists 
$\beta = \beta(b, \ell, \theta_{\min},\alpha) > 0$ such that the
following hold.
Let $T$ be an $\ell$-level balanced $b$-ary tree, and consider the 
$CFN(\theta,\eta)$ model on $T$, where $\min \theta \geq \theta_{\min}$
and $\min \eta \geq \eta_{\min}$. Then
\begin{equation} \label{eq:maj_prop}
\widehat{\Maj}(\theta,\eta) \geq \min\{\alpha \eta_{\min}, \beta\}.
\end{equation}

In particular, given $b$ and $\theta_{\min}$ such that
$b \theta_{\min}^2 > g^2 > 0$, there exist $\ell(b,\theta_{\min})$,
$\alpha(b,\theta_{\min}) > g^{\ell}$ and $\beta(b,\theta_{\min}) > 0$, such
that any $CFN(\theta,\eta)$ model on the $\ell$-level $b$-ary tree
satisfying $\min \theta \geq \theta_{\min}$ and $\min \eta \geq
\eta_{min}$ must also satisfy (\ref{eq:maj_prop})
\end{Theorem}

Theorem \ref{thm:maj_good} is a generalization of Lemma 
\ref{lem:recur_fix_eps}. 

\prooft{of Lemma \ref{lem:recur_fix_eps} \cite{M:recur}}
Follows immediately from the second assertion in Theorem
\ref{thm:maj_good}, where $g=1$, and $\eta_0$
(of Lemma \ref{lem:recur_fix_eps}) is chosen between $0$ and $\beta$ (of Theorem \ref{thm:maj_good}).
$\QED$

The following lemma is a generalization of the second claim in
Lemma \ref{lem:rec_fix_eps}.

\begin{Lemma} \label{cor:rec_maj}
Let $b$ and $\theta_{\min}$ be such that $b \theta_{\min}^2 > g^2 > 0$.
Let $\ell(b,\theta_{\min})$,
$\alpha(b,\theta_{\min}) > g^{\ell}$ and 
$\beta(b,\theta_{\min}) > 0$, be such
that (\ref{eq:maj_prop}) holds.

Consider the $CFN(\theta,\eta)$ model on the family of 
balanced tree of $q$ levels, where all
internal nodes have at least $b$ children, a total of $n$ leaves,  
$\min \theta \geq \theta_{\min}$ and $\min \eta \geq \beta$. 

Then for all $T$ and $i \geq 0$, there exists a map 
$\Psi = \Psi_T  
: \{\pm 1\}^{\partial T} \to \{\pm 1\}^{L_{\partial - i \ell}}$ for
which the following hold.

If $\sigma$ is distributed according to the $CFN(\theta,\eta)$
model on $T$, and $\sigma' = \Psi(\sigma_{\partial T})$, then
$(\sigma'_v)_{v \in L_{\partial - i \ell}} = 
(\sigma_v \tau_v)_{v \in L_{\partial - i \ell}}$, where
$\tau_v$ are independent variables.
Moreover, $\tau_v$ are  
independent of $(\sigma_v)_{d(v,\partial T) \geq i \ell}$ and satisfy
$\E[\tau_v] \geq \min\{1,g^{i \ell}\} \beta$ for all $v$.

The map $\Psi$ may be constructed from the $(i \ell)$-topology of $T$.
In particular, if $T_1$ and $T_2$ have the same $(i \ell)$-topology, then 
$\Psi_{T_1} = \Psi_{T_2}$. Furthermore, $\Psi(\sigma_{\partial T})$ 
is computable in time polynomial in $n$.

 
\end{Lemma}

\proofs
Similarly to Lemma \ref{lem:rec_fix_eps},
the proof is by induction on $i$. 
The claim is trivial for $i = 0$.

For the induction step, suppose that we are given $d^{\ast}_{i \ell+\ell}$.
Label all the vertices $v \in \cup_{j=0}^{(i+1) \ell} L_{\partial - j}$ by
the $(i+1) \ell$-labeling of $T$. 

By the induction hypothesis, there exists a map $\Psi'$ such that 
$\Psi'(\sigma_{\partial T}) = 
(\sigma_v \tau_v)_{v \in L_{\partial T - i \ell}}$, 
where $\tau_v$ are independent ${\pm 1}$ variables, 
 independent of $(\sigma_v)_{d(v,\partial T) \geq i \ell}$ and satisfy
 $\E[\tau_v] \geq  \min\{1,g^{i \ell}\} \beta$ for all $v$.  

By the properties of the labeling, for each 
$w \in L_{\partial - (i+1)\ell}$, we can find a set 
$R(w) \subset L_{\partial - i \ell}$, 
which is the set of leaves of an
$\ell$-level $b$-ary tree rooted at $w$.

We now let  $\Psi(\sigma_{\partial T}) = 
(\widehat{\sigma}_w)_{w \in L_{\partial - i \ell - \ell}}$, where
\[
\widehat{\sigma}_w = 
\Maj((\Psi'(\sigma_{\partial T}))_v : v \in R(w)) = 
\Maj(\sigma_v \tau_v : v \in R(w)).
\]
By Theorem \ref{thm:maj_good} it follows that 
$(\widehat{\sigma}_w)_{w \in L_{\partial - i \ell - \ell}} =  
(\sigma_w \tau_w)_{w \in L_{\partial - i \ell - \ell}}$,
where $\tau_w$ are independent $\pm 1$ variables. 
Moreover, $\tau_w$ are independent 
of $(\sigma_v)_{d(v,\partial T) \geq i \ell + \ell}$ and 
satisfy $\E[\tau_w] \geq  \min\{1,g^{i \ell + \ell}\} \beta$, 
for all $w$, proving the second
claim.

Note that in order to compute the function $\Psi$,
one applies the majority function recursively starting at a subset of
the leaves. Therefore $\Psi$
is computable in time polynomial in $n$.
$\QED$

The following Lemma shows why the second 
assertion of the Theorem \ref{thm:maj_good} follows from the first one.
\begin{Lemma} \label{lem:stirling}
\[
\lim_{\ell \to \infty} \frac{a(b^{\ell}) \theta^{\ell}}
                {\sqrt{\frac{2}{\pi}} b^{\ell/2} \theta^{\ell}} = 1.
\]
In particular, if $b \theta^2 > g^2$, then 
$a(b^{\ell}) \theta^{\ell} > g^{\ell}$, for all sufficiently large $\ell$.
\end{Lemma}
\proofs
Stirling's formula implies that
\[
a(d) = 2^{1-d} \lceil \frac{d}{2} \rceil \binom{d}{\lceil \frac{d}{2} \rceil}
= (1 + o(1)) \sqrt{\frac{2}{\pi}} \sqrt{d}.
\]
Now the claim follows.
$\QED$


The role that $a(d)$ plays for the majority algorithm
 is presented in the following lemma.

\begin{Lemma} \label{lem:bin_prob}
\begin{enumerate}
\item
Let $X,Y_1,\ldots,Y_d$ be a sequence of $\pm 1$ random variables such that $Y_2,\ldots,Y_d$ are i.i.d. with
$\E[Y_i] = 0$, $\E[Y_1 | X = 1] = -\E[Y_1 | X = -1] = \theta$, and $Y_2,\ldots,Y_d$ are independent of $X,Y_1$.
Then
\begin{equation} \label{eq:easy_conc}
\E[X \Maj(Y_1,\ldots,Y_d) | X = 1] =  \E[X \Maj(Y_1,\ldots,Y_d) | X = -1] = \theta \frac{a(d)}{d}.
\end{equation}
where $a(d)$ is given in (\ref{eq:def_a}).
\item
Let $X,Y_1,\ldots,Y_{d-1}$ be a collection of random variables, 
where $X$ is non-negative, $Y_1,\ldots,Y_{d-1}$ are symmetric and
\begin{itemize}
\item
$Y_1,\ldots,Y_{d-1}$ are independent, and
\item
$\P[X \geq \max_i |Y_i|] = 1$.
\end{itemize}
Then
\begin{equation} \label{eq:hard_conc}
\E[\sign(X + \sum_{i=1}^{d-1} Y_i)] \geq \frac{a(d)}{d}.
\end{equation}
\end{enumerate}
\end{Lemma}
\proofs
\begin{enumerate}
\item
Let $\tilde{Y}$ be a $\pm 1$ variable which is independent of 
$X,Y_2,\ldots,Y_d$, and $\E[\tilde{Y}] = 0$. Let $Z$ be a random
variable which is independent of $X,Y_2,\ldots,Y_d,\tilde{Y}$,
such that $\P[Z = 1] = \theta$, and $\P[Z = 0] = 1 - \theta$.
Note that $(Y_1,\ldots,Y_d)$ and 
$(Z X + (1 - Z) \tilde{Y},Y_2,\ldots,Y_d)$ have the same distribution.
Therefore,
\begin{eqnarray} \label{eq:decY}
\E[X \Maj(Y_1,\ldots,Y_d) | X = 1] &=& 
\theta \E[X \Maj(X,Y_2,\ldots,Y_d) | X = 1] + \\ \nonumber &\,& 
(1 - \theta) \E[X \Maj(\tilde{Y},Y_2\ldots,Y_d) | X = 1].
\end{eqnarray}
Since $\tilde{Y},Y_2,\ldots,Y_d$ are independent of $X$, it follows that
\[
\E[X \Maj(\tilde{Y},Y_2,\ldots,Y_d) | X = 1] = 0.
\]
Therefore by (\ref{eq:decY}), in order to prove the Lemma, it suffices to
show that if $Y_2,\ldots,Y_d$
are i.i.d. $\pm 1$ random variables with $\E[Y_i] = 0$, then
\begin{equation} \label{eq:majY}
\E[\Maj(1,Y_2,\ldots,Y_d)]  =
\frac{a(d)}{d}.
\end{equation}
It is helpful to note that
$\E[\Maj(1,Y_2,\ldots,Y_d)]  = \P[\sign(1 + \sum_{i=2}^d Y_i) \neq  \sign(\sum_{i=2}^d Y_i)]$.
There are two cases to consider:
\begin{itemize}
\item
$d=2e+1$ is odd:
\[
\P[\sign(1 + \sum_{i=2}^d Y_i) \neq  \sign(\sum_{i=2}^d Y_i)] =
\P[\sum_{i=2}^d Y_i = 0] = 2^{-2e} \binom{2e}{e} = \frac{a(d)}{d}.
\]
\item
$d=2e$ is even:
\[
\P[\sign(1 + \sum_{i=2}^d Y_i) \neq  \sign(\sum_{i=2}^d Y_i)] =
\P[\sum_{i=2}^d Y_i = -1]
= 2^{-2e+1} \binom{2 e - 1}{e} = \frac{a(d)}{d}.
\]
\end{itemize}
\item
The general case follows by conditioning 
from the case where for all $i$, $Y_i$ is a $\pm y_i$ random variable,
$X = x$ is a constant and $x \geq \max_{1 \leq i \leq d-1} |y_i|$.

If $x = |y_1| = \ldots = |y_{d-1}|$ then the claim follows from the
proof of the first part of the lemma. We may therefore assume the 
strict inequality $x > |y_{d-1}|$.

We now show that it suffices to prove the claim for 
$x,y_1,\ldots,y_{d-1}$ such that
\begin{equation} \label{eq:no_ties} 
\P[x + \sum_{i=1}^{d-1} Y_i = 0] = 0.
\end{equation}
Indeed, if $\P[x + \sum_{i=1}^{d-1} Y_i = 0] > 0$, let 
\[
\eps = \frac{1}{2} \min \{ |c_0 x + \sum_{i=1}^{d-1} c_i y_i| :
             c_0 x + \sum_{i=1}^{d-1} c_i y_i \neq 0 
             \mbox{ and } c_i \in \{-1,0,1\} \mbox{ for all } i\}.
\]
Define $y^{+}_i,y^{-}_i$ for $1 \leq i \leq d-1$ by
\[
y^{\pm}_i = \left\{\begin{array}{ll}
                  y_i & \mbox{if }i < d-1, \\ 
                  y_i \pm \eps & \mbox{if } i = d-1.
                  \end{array}
                  \right.
\]
Let $Y^{+}_i$ be independent symmetric $\pm y^{+}_i$ variables and 
let $Y^{-}_i$ be independent symmetric $\pm y^{-}_i$ variables.
Note that
\begin{equation} \label{eq:tie_breaker}
\E[\sign(x + \sum_{i=1}^{d-1} Y_i)] = 
\frac{1}{2} \E[\sign(x + \sum_{i=1}^{d-1} Y^{+}_i)] +
\frac{1}{2} \E[\sign(x + \sum_{i=1}^{d-1} Y^{-}_i)],
\end{equation}
$x \geq \max_{1 \leq i \leq d-1} |y^{\pm}_i|$, $x >
|y^{\pm}_{d-1}|$ and 
\begin{equation} \label{eq:dead_ties}
\P[x + \sum_{i=1}^{d-1} Y^{+}_i = 0] =
\P[x + \sum_{i=1}^{d-1} Y^{-}_i = 0] = 0.
\end{equation}
From (\ref{eq:tie_breaker}) and (\ref{eq:dead_ties}) it follows that 
we may assume (\ref{eq:no_ties}).

By (\ref{eq:no_ties}) and symmetry
\begin{equation} \label{eq:symiso}
\E[\sign(x + \sum_{i=1}^{d-1} Y_i)] = 
1 - 2 \P[\sum_{i=1}^{d-1} Y_i < -x] = 
1 -  \P[\sum_{i=1}^{d-1} Y_i < -x] - \P[\sum_{i=1}^{d-1} Y_i > x]. 
\end{equation}
Let $U_{-} \subset \{-1,1\}^{d-1}$ be defined as 
$U_{-} = \{(b_1,\ldots,b_{d-1}) : \sum b_i y_i < -x\}$
and $U_{+} = \{(b_1,\ldots,b_{d-1}) : \sum b_i y_i > x\}$.
Rewriting (\ref{eq:symiso}), we get 
\begin{equation} \label{eq:symiso2}
\E[\sign(x + \sum_{i=1}^{d-1} Y_i)] = 1 - \P[U_{-}] - \P[U_{+}], 
\end{equation}
where $\P$ is the uniform measure on $\{-1,1\}^{d-1}$.
Note that if $h$ denotes the Hamming distance, then
\[
h(U_{+},U_{-}) := \min_{b_{+} \in U_{+}, b_{-} \in U_{-}} h(b_{+},b_{-})
\geq 2.
\]
By the isoperimetric inequality for the discrete cube \cite{FF} 
(see \cite{Harper}; \cite{Bezrukov} for background) 
it follows that maximizers of the following quantity 
\[
\max\{ p : \P[U] = \P[U'] = p; U,U' \subset \{-1,1\}^{d-1} \mbox{ and } 
h(U,U') \geq 2\}
\]
are obtained as follows.
\begin{itemize}
\item
If $d=2e+1$ is odd, then the maximum is obtained for 
\[
U  = \{b : \sum_{i=1}^{d-1} b_i \geq 2\},
\]
\[
U' = \{b : \sum_{i=1}^{d-1} b_i \leq -2\}.
\]
Therefore
\begin{equation} \label{eq:iso1}
1 - \P[U_{+}] - \P[U_{-}] \geq 
1 - \P[U] - \P[U'] = 2^{-2e} \binom{2e}{e} =
\frac{a(d)}{d}.
\end{equation}
\item
If $d=2e$ is even, then the maximum is obtained for
\[
U = \{b : \sum_{i=1}^{d-1} b_i \geq 3\} \cup
     \{b : \sum_{i=1}^{d-1} b_i = 1 \mbox{ and } b_1 = b_2 = 1\},
\]
\[
U' = \{b : \sum_{i=1}^{d-1} b_i \leq -3\} \cup
     \{b : \sum_{i=1}^{d-1} b_i = -1 \mbox{ and } b_1 = b_2 = -1\}.
\]
Therefore in this case 
\begin{equation}
\label{eq:iso2}
1 - \P[U] - \P[U'] = 2 \times 2^{-2e+1} 
\left(\binom{2e-1}{e} - \binom{2e-3}{e - 1} \right) \geq 
2^{-2e+1} \binom{2e-1}{e} = \frac{a(d)}{d}. 
\end{equation}
(the assumption $|y_{d-1}| < x$ resulted in a better bound here than in
$\pm 1$ case).
\end{itemize}
Now (\ref{eq:symiso2}) follows from (\ref{eq:iso1}) and 
(\ref{eq:iso2}).
\end{enumerate}
$\QED$


\begin{Lemma} \label{lem:maj_crit}
Let $T$ be the $\ell$-level $b$-ary tree.
Suppose that $\alpha < a(b^{\ell}) \theta_{\min}^{\ell}$.
Then there exists $\eps = \eps(b,\ell,\theta_{\min},\alpha) > 0$, 
s.t. if $\max \eta \leq \eps$, 
and $\min \theta \geq \theta_{\min}$, then the $CFN(\theta,\eta)$ model
on $T$ satisfies
\begin{equation} \label{eq:crit_f}
\widehat{\Maj}(\theta,\eta)  \geq
\alpha \frac{\sum_{v \in \partial T} \eta(v)}{b^{\ell}}.
\end{equation}
\end{Lemma}

\proofs
In order to prove (\ref{eq:crit_f}), it suffices to show that
for $\alpha < a(b^{\ell}) \theta_{\min}^{\ell}$, there exists $\eps >
0$, s.t. if $\max \eta \leq \eps$,
then for all $v$
\begin{equation} \label{eq:der_f}
\frac{\partial \widehat{\Maj}}{\partial \eta(v)}(\theta,\eta) \geq
\frac{\alpha}{b^{\ell}}.
\end{equation}
Note that $\widehat{\Maj}(\theta,\eta)$ is a polynomial in $\theta$ and $\eta$.
Therefore all the derivatives of $\widehat{\Maj}(\theta,\eta)$
with respect to the $\theta$ and $\eta$ variables are
uniformly bounded.
In particular, there exists a constant $C$ (which depends on $b$ and
$\ell$ only) such that for all $\theta$
and $\eta$ we have that
\begin{equation} \label{eq:comp_eta_to_eta0}
|\frac{\partial \widehat{\Maj}}{\partial \eta(v)}(\theta,\eta) - 
 \frac{\partial \widehat{\Maj}}{\partial \eta(v)}(\theta,0)| \leq 
C \max_i \eta_i.
\end{equation}
Therefore it suffices to show that for all
$\theta$ with $\min_e \theta(e) \geq \theta_{\min}$, we have that
\begin{equation} \label{eq:der_0}
\frac{\partial \widehat{\Maj}}{\partial \eta(v)}(\theta,0)
\geq \frac{a(b^{\ell}) \theta_{\min}^{\ell}}{b^{\ell}}.
\end{equation}
By (\ref{eq:comp_eta_to_eta0}) this will imply (\ref{eq:der_f}) for all $\eta$
satisfying $\max_v \eta(v) \leq \eps$, where 
\[
\eps = \frac{a(b^{\ell}) \theta_{\min}^{\ell} - \alpha}{C b^{\ell}}.
\]
Fix $v \in \partial T$, and let $e_1,\ldots,e_{\ell}$ 
be the path in $T$ from the root $\rho$ to $v$.
Let $\gamma \in [0,1]$ and
\[
\eta_0(w) = \left\{\begin{array}{ll}
                   0 & \mbox{ if } w \neq v, \\
                   \gamma & \mbox{ if } w = v.
                   \end{array} \right.
\]
$\widehat{\Maj}(\theta,\eta_0)$ is the covariance of the majority of $b^{\ell}$
i.i.d. $\pm 1$ variables and $\sigma_{\rho}$. 
One of these variables, $\sigma_v$, satisfies
$\E[\sigma_v \sigma_{\rho}] = \gamma \prod_{i=1}^{\ell} \theta(e_i)$, 
while all
the other variables are independent of $\sigma_{\rho}$ and $\sigma_v$.
Therefore by part 1 of Lemma (\ref{lem:bin_prob}) it follows that
\[
\widehat{\Maj}(\theta,\eta_0) =
\frac{a(b^{\ell})}{b^{\ell}} \gamma \prod_{i=1}^{\ell} \theta(e_i). 
\]
So
\[
\frac{\partial \widehat{\Maj}}{\partial \eta(v)}(\theta,0) =
\frac{a(b^{\ell})}{b^{\ell}} \prod_{i=1}^{\ell} \theta(e_i). 
\]
The assumption on $\theta$ implies that $\prod_{i=1}^{\ell} \theta(e_i) \geq \theta_{\min}^{\ell}$.
Therefore,
\[
\frac{\partial \widehat{\Maj}}{\partial \eta(v)}(\theta,0) \geq
\frac{a(b^{\ell}) \theta_{\min}^{\ell}}{b^{\ell}},
\]
to obtain (\ref{eq:der_0}), as needed.
$\QED$

\begin{Lemma} \label{lem:maj_far}
Let $T$ be an $\ell$-level balanced $b$-ary tree.
Suppose that $\min \theta \geq \theta_{\min}$, and $\max \eta \geq \eta_{\max}$,
then
\begin{equation} \label{eq:non_crit_f}
\widehat{\Maj}(\theta,\eta) \geq \frac{a(b^{\ell})}{2^{\ell} b^{2 \ell + 1}}
h(\theta_{\min})^{\ell-1} h(\theta_{\min} \eta_{\max}),
\end{equation}
where
\begin{equation} \label{eq:def_h}
h(x) = \min\{1,\frac{x}{1-x}\}.
\end{equation}
\end{Lemma}

\proofs
We use the ``random cluster'' representation of the model 
(see \cite{Grimmett} for background on
percolation and random-cluster models).
Declare an edge $e=(u,v)$ {\em open} with probability $\theta(e)$ if $e$
is not adjacent to
$\partial T$, and with probability $\theta(e) \eta(v)$ if $v \in \partial T$, independently
for all edges.
An edge which is not open is declared {\em closed}. Given the clusters (connected open components)
of the random cluster representation, color the root cluster by the root color, and
each of the other clusters, by an independent unbiased $\pm 1$ variable.
This gives the same distribution on coloring as the original coloring procedure
(this should be clear; see e.g. \cite{M:lam2} for more details).

Assume that the root color is $1$, and let ${\C}_{\rho}$ be the root cluster. Let
\[
X = \sum_{v \in {{\C}'}_{\rho}} \sigma_v = |{{\C}'}_{\rho}| = |{\C}_{\rho} \cap \partial T|,
\]
where ${{\C}'}_{\rho} = {\C}_{\rho} \cap \partial T$.
Let ${\C}_1,\ldots,{\C}_K$ be all other clusters (note that $K$ is a random variable) and let
\[
Y_i = \sum_{v \in {\C}_i'} \sigma_v,
\]
where ${\C}_i' = {\C}_i \cap \partial T$.
Conditioned on ${\C}_{\rho},{\C}_1,\ldots,{\C}_K$,
\[
\P[Y_i = \pm |\C_i \cap \partial T|] = 1/2,
\]
and the $Y_i$'s are independent conditioned on 
${\C}_{\rho},{\C}_1,\ldots,{\C}_K$.
Clearly,
\begin{equation} \label{eq:f_as_sum}
\widehat{\Maj}(\theta,\eta) = \E[\sign(X + \sum_{i=1}^K Y_i)].
\end{equation}

Note that conditioned on ${\C}_{\rho},{\C}_1,\ldots,{\C}_K$, the
variable $\sum_{i = 1}^K Y_i$ is symmetric and therefore,
\begin{equation} \label{eq:bad_cond_good}
\E[\sign(X + \sum_{i=1}^K Y_i) \,\, | \,\,|{{\C}'}_{\rho}| < \max_i |{{\C}'}_i|] \geq 0.
\end{equation}

Moreover, below we prove that
\begin{equation} \label{lem:c_is_big}
\P[|{{\C}'}_{\rho}| \geq \max_i |{{\C}'}_i|] \geq 2^{- \ell}
b^{-\ell-1} h(\theta_{\min})^{\ell-1} h(\theta_{\min} \eta_{\max}).
\end{equation}
When $X > 0$, there are at most $b^{\ell} - 1$ non-zero variables among the $Y_i's$.
Therefore part 2 of Lemma \ref{lem:bin_prob} implies that
\begin{equation} \label{eq:sum_is_big}
\E[\sign(X + \sum_{i=1}^K Y_i) \,\, | \,\,|{{\C}'}_{\rho}| \geq \max_i |{{\C}'}_i|| \geq
\frac{a(b^{\ell})}{b^{\ell}}.
\end{equation}
Combining (\ref{eq:sum_is_big}) and (\ref{eq:bad_cond_good}) via 
(\ref{eq:f_as_sum}) we obtain:
\begin{eqnarray*}
\widehat{\Maj}(\theta,\eta) &=&
\P[|{{\C}'}_{\rho}| < \max_i |{{\C}'}_i|]\,
\E[\sign(X + \sum_{i=1}^K Y_i) \,\, |\,\, |{{\C}'}_{\rho}| < 
   \max_i |{{\C}'}_i|] \\ &+&
\P[|{{\C}'}_{\rho}| \geq \max_i |{{\C}'}_i|]\,\,
\E[\sign(X + \sum_{i=1}^K Y_i) \, |\,\, |{{\C}'}_{\rho}| 
\geq \max_i |{{\C}'}_i|] \\ &\geq&
\frac{a(b^{\ell})}{2^{\ell} b^{2 \ell + 1}} h(\theta_{\min})^{\ell-1} h(\theta_{\min} \eta_{\max}),
\end{eqnarray*}
as needed.

It remains to prove (\ref{lem:c_is_big}).
Let $\Omega$ be the probability space of all random cluster configurations.
We prove (\ref{lem:c_is_big}) by constructing a map $G: \Omega \to \Omega$ such that for all $\omega \in \Omega$,
\begin{equation} \label{eq:Glarge}
|{{\C}'}_{\rho}(G(\omega))| \geq \max_i |{{\C}'}_i(G(\omega))|,
\end{equation}
and for all $\omega$
\begin{equation} \label{eq:G_RN}
\P[G^{-1}(\omega)] \leq b^{\ell + 1}\,
2^{\ell}\,h(\theta_{\min})^{1-\ell}\,
h(\theta_{\min} \eta_{\max})^{-1}\, \P[\omega].
\end{equation}

If $\omega$ satisfies $|{{\C}'}_{\rho}| \geq \max_i |{{\C}'}_i|$, then we let $G(\omega) = \omega$.
Otherwise, let $\C$ be a cluster such $|\C \cap \partial T| = \max_i |{{\C}'}_i|$.
If $\max_i |{{\C}'}_i| = 1$, we let $\C$ be a cluster which contains a $v \in \partial T$ with $\eta(v) \geq \eta_{\max}$.
Let $u \in \C$ be the vertex closest to the root $\rho$.
Let $G(\omega)$ be the configuration which is obtained from $\omega$ be setting all the edges on the
path from $u$ to $\rho$ to be open.

It is clear that $G(\omega)$ satisfies (\ref{eq:Glarge}).
Let $\omega$ be such that $|{{\C}'}_{\rho}(\omega)| \geq \max_i |C'_i(\omega)|$.
Then any element in $G^{-1}(\omega)$ is obtained by
\begin{itemize}
\item
choosing a vertex $u \in T$ such that either $u \notin \partial T$ or $u \in \partial T$ and
$\eta(u) \geq \eta_{\max}$.
\item
choosing a subset $S$ of the edges on the path from $u$ to $\rho$
\item
Setting all the edges of $S$ to be close.
\end{itemize}
If $\omega'$ is the configuration thus obtained, then clearly,
\begin{equation} \label{eq:inv_clust}
\P[\omega'] \leq h(\theta_{\min})^{1-\ell} h(\theta_{\min}
\eta_{\max})^{-1} \P[\omega].
\end{equation}
It remains to count the number of $\omega'$ which may be obtained from 
$\omega$.
There are at most $b^{\ell+1}$ choices for $u$.
Moreover, there are at most $2^{\ell}$ subsets of the edges we want to update at the second stage.
Thus there are at most $b^{\ell+1} 2^{\ell}$ pre-images $\omega'$ to consider, each satisfying (\ref{eq:inv_clust}).
We thus obtain (\ref{eq:G_RN}) as needed.
$\QED$

\prooft{of Theorem \ref{thm:maj_good}}
By Lemma \ref{lem:maj_crit}, there exist $\eps > 0$ 
such that if for all $v \in \partial T$,
it holds that $\eta_{\min} \leq \eta(v) \leq \eps$, then
\[
\widehat{\Maj}(\theta,\eta) \geq \alpha \eta_{\min}.
\]
By Lemma \ref{lem:maj_far}, it follows that if 
$\max_{v \in \partial T} \eta(v) \geq \eps$, then
\[
\widehat{\Maj}(\theta,\eta) \geq \beta,
\]
where
\[
\beta = \frac{a(b^{\ell})}{2^{\ell} b^{2 \ell + 1}}
h(\theta_{\min})^{\ell-1} 
h(\theta_{\min} \eps),
\]
and $h$ is given by (\ref{eq:def_h}).
Now the first claim follows.
The second claim follows from the first claim by Lemma \ref{lem:stirling}.
$\QED$

\section{Four point condition and topology} \label{sec:four_points}

In this section we discuss how to reconstruct 
the $\ell$-topology of a balanced tree, given
the correlation between colors at different leaves. 
The analysis in this section does not exhibit
a phase transition when $b \theta_{\min}^2 = 1$, 
as $\theta_{\min}$ has a continuous role
in the bounds below. 
We follow the well known technique of ``$4$-point condition''. 
However, as we require
to reconstruct only the local topology, and consider only balanced trees,
the number of samples needed is {\em logarithmic} in $n$.

The following theorem generalizes the first part of 
Lemma \ref{lem:rec_fix_eps}.

\begin{Theorem} \label{thm:4point}
Let $\ell$ be a positive integer and 
consider the $CFN(\theta,\eta)$ model on the family of 
balanced tree on $n$ leaves, where
\begin{itemize}
\item[I.]
For all edges $e$, 
$\theta_{\min} \leq \theta(e) \leq \theta_{\max}$, where $\theta_{\min} > 0$
and $\theta_{\max} < 1$, and
\item[II.]
For all $v \in \partial T$, $\eta_{\min} \leq \eta(v)$,
where $\eta_{\min} > 0$.
\end{itemize}
Then there exists a map $\Phi$ from $\{\pm 1\}^{k n}$ to 
the space of $\ell$-topologies on $n$ leaves, 
such that,
\[
\P[\Phi \left((\sigma^t_{\partial T})_{t=1}^k \right) 
= \ell-\mbox{topology of } T]
\geq 1 - \delta,
\]  
where $(\sigma^t_{\partial T})_{t=1}^k$  
are $k$ independent samples of the process at the leaves of $T$,
and 
\begin{equation} \label{eq:4_point_error}
\delta \leq 
n^2 \exp(- c^{\ast}\, k \, \theta_{\min}^{8 \ell + 8} \eta_{\min}^8 
(1 - \theta_{\max})^2),
\end{equation}
with $c^{\ast} \geq 1/2048$.
Moreover, $\Phi$ is computable in polynomial time in $n$ and $k$.
\end{Theorem}

\begin{Definition} \label{def:D}
Consider the $CFN(\theta,\eta)$ coloring of a tree $T$.
For any two leaves $u,v$, let
\[
\theta(u,v) = \eta(u) \eta(v) \prod_{w \in \pa(u,v)} \theta(w),
\]
and
\[
D(u,v) = -\log \theta(u,v) = -\log (\eta(u)) - \log(\eta(v))
- \sum_{w \in \pa(u,v)} \log(\theta(w)).
\]
\end{Definition}

\begin{Lemma} \label{lem:d_D}
Suppose we are given $k$ samples $(\sigma^t_{\partial T})_{t=1}^k$, of the
$CFN(\theta,\eta)$ coloring of a tree $T$ as in Theorem \ref{thm:4point}.
For $u,v \in \partial T$, let
\[
c(u,v) = \frac{1}{k}\sum_{t=1}^k \sigma^t_u \sigma^t_v,
\]
(note that the expected value of $c(u,v)$ is $\theta(u,v)$),
and
\[
D^{\ast}(u,v) = \left\{ \begin{array}{ll} 
                -\log c(u,v) & \mbox{if } c(u,v) > 0, \\
                \infty       & \mbox{if } c(u,v) \leq 0.
                \end{array} 
                \right.
\]
Let $\theta_{\ast} = \theta_{\min}^{2 \ell + 2} \eta_{\min}^2/2$. 
Define $uRv$ if $c(u,v) \geq \theta_{\ast}$ and $uR'v$ if 
$c(u,v) \geq \frac{15}{8} \theta_{\ast}$. 
($uRv$ roughly means that $u$ and $v$ are close;
$uR'v$ means that $u$ and $v$ are even closer) and let
$1/4 \geq \eps > 0$. Then with probability at least 
\begin{equation} \label{eq:error_eps}
1 - n^2  \exp(-k \theta_{\ast}^4 \eps^2/8),
\end{equation}
\begin{itemize}
\item[I.]
$uR'v$ for all $u$ and $v$ such that $d(u,v) \leq 2 \ell +
2$. 
\item[II.]
For all $u$ and $v$, 
if there exists a $w$ such that $uRw$ and $vRw$ then
$|D(u,v) - D^{\ast}(u,v)| < \eps$.
\end{itemize}
\end{Lemma}
\proofs
Define
\[
A = \{ \exists (u,v) \mbox{ s.t. } 
|c(u,v) - \theta(u,v)| \geq \alpha \},
\]
where $\alpha = \eps \theta_{\ast}^2/2$.
We claim that conditioned on $A^c$, both I. and II. hold.
If $u,v  \in \partial T$ satisfy $d(u,v) \leq 2 \ell + 2$, then
$\theta(u,v) \geq 2 \theta_{\ast}$.
Therefore conditioned on $A^c$, all 
$u,v \in \partial T$ s.t. $d(u,v) \leq 2 \ell + 2$, 
must satisfy  
\[
c(u,v) \geq 2 \theta_{\ast} - \eps \theta_{\ast}^2/2 \geq
\frac{15}{8} \theta_{\ast},
\]
and I. follows.

Conditioned on $A^c$, if $uRw$, then 
$c(u,w) \geq \theta_{\ast}$ and therefore $\theta(u,w) > \theta_{\ast} - \alpha$. 
Similarly, if $vRw$, then $\theta(v,w) > \theta_{\ast} - \alpha$.
Now
\begin{eqnarray*}
\theta(u,v) &=& 
\eta(u) \eta(v) \prod_{y \in \pa(u,v)} \theta(y) \geq
\left( \eta(u) \eta(w) \prod_{y \in \pa(u,w)} \theta(y) \right)
\left( \eta(w) \eta(v) \prod_{y \in \pa(w,v)} \theta(y) \right) \\ &=& 
\theta(u,w) \theta(w,v)
> (\theta_{\ast} - \alpha)^2,
\end{eqnarray*}
and therefore conditioned on $A^c$,
\[
c(u,v) > (\theta_{\ast} - \alpha)^2 - \alpha. 
\]
Therefore conditioned on $A^c$, by the mean value theorem,
\begin{eqnarray*}
|D^{\ast}(u,v) - D(u,v)| &=& |\log c(u,v) - \log \theta(u,v)| \leq
\frac{|c(u,v) - \theta(u,v)|}{(\theta_{\ast} - \alpha)^2 - \alpha} <  
\frac{\alpha}{(\theta_{\ast} - \alpha)^2 - \alpha} 
\\ \nonumber
&=&
\frac{\eps \theta_{\ast}^2/2}
{(\theta_{\ast} - \eps \theta_{\ast}^2/2)^2 - \eps \theta_{\ast}^2/2} 
= \frac{\eps}{2(1 - \eps \theta_{\ast}/2)^2 - \eps} < \eps
\end{eqnarray*}
to obtain II.

By Lemma  \ref{lem:dev_fix_eps},
$\P[A]$ is bounded by 
\begin{equation} \label{eq:longer_dev}
\P[A] \leq \binom{n}{2} 2 \exp(-k \theta_{\ast}^4 \eps^2 / 8) \leq
n^2  \exp(-k \theta_{\ast}^4 \eps^2 / 8),
\end{equation}
as needed.

\QED

For a set $V$ of size $4$ a {\em split} is defined as a partition of
$V$ into two sets of size $2$. We will write $v_1 v_2 | v_3 v_4$ 
for the split $\{ \{v_1,v_2\}, \{v_3, v_4\} \}$. Note that a
$4$-element set has exactly $3$ different splits.

\begin{Lemma} \label{lem:top}
Let $T=(V,E)$ be a balanced tree. Let 
$\Delta : E \to \R_{+}$ be a positive function.
For $u,v \in V$, let 
$\Delta(u,v) = \sum_{e \in \path(u,v)} \Delta(e)$. 
For a split $\Gamma = u_1 u_2 | u_3 u_4$, let
\[
\Delta(\Gamma) = \Delta(u_1,u_2) + \Delta(u_3,u_4).
\]
Then
\begin{itemize}
\item
If $\Gamma_1$ and $\Gamma_2$ are two splits of
$\{u_1,u_2,u_3,u_4\}$, then either $\Delta(\Gamma_1) = \Delta(\Gamma_2)$, or
$|\Delta(\Gamma_1) - \Delta(\Gamma_2)| \geq 2 \Delta_{\min}$, where
\[
\Delta_{\min} = \min \{\Delta(e) : e \mbox{ not adjacent to }
\partial T\}.
\]
\item
Let $R$ be a binary relation on $\partial T$
such that $uRv$ whenever $d(u,v) \leq 2 \ell + 2$. Write $R(v)$
for the set of elements which are related to $v$.
Then in order to reconstruct the $\ell$-topology of the tree
it suffices to find for all $u \in \partial T$ 
and all $\{u_1,u_2,u_3,u_4\} \subset R(u)$,
all minimizers of 
\[
\{ \Delta(\Gamma) :  \Gamma \mbox{ a split of } \{u_1,u_2,u_3,u_4\} \}
\]
(we call such minimizers minimal splits).
\end{itemize}
\end{Lemma}

\proofs
Let $U$ be a set of four vertices. Note that either there is a unique
split $u_1 u_2 | u_3 u_4$ of $U$ such that 
$\pa(u_1,u_2) \cap \pa(u_3,u_4)$ is empty, or for all splits 
$u_1 u_2 | u_3 u_4$, the set $\pa(u_1,u_2) \cap \pa(u_3,u_4)$ consists
of a single vertex.

Suppose that $u_i \in \partial T$, for $1 \leq i \leq 4$
and that $\pa(u_1,u_2) \cap \pa(u_3,u_4)$ is empty.
Let $u_{1,2}$ be the point on $\pa(u_1,u_2)$ which is closest
to $\pa(u_3,u_4)$. Define $u_{3,4}$ similarly.
Then
\begin{equation} \label{eq:for1}
\Delta(u_1,u_3) + \Delta(u_2,u_4) = 
\Delta(u_1,u_4) + \Delta(u_2,u_3),
\end{equation}
and
\begin{equation} \label{eq:for2}
\Delta(u_1,u_3) + \Delta(u_2,u_4) - 
\Delta(u_1,u_2) - \Delta(u_3,u_4) =
2 \sum_{e \in \pa(u_{1,2},u_{3,4})} \Delta(e) \geq 2 \Delta_{\min}.
\end{equation}
If on the other hand, $\pa(u_1,u_2) \cap \pa(u_3,u_4)$ consists of
a single point, then for all permutations ${i,j},{k,\ell}$ of $1,2,3,4$,
\begin{equation} \label{eq:for3}
\Delta(u_i,u_j) + \Delta(u_k,u_{\ell}) \mbox{ has the same value}.
\end{equation}
The first claim follows.

Let $\rho$ be the root of the tree and let $q$ be the distance between 
$\rho$ and the leaves (since the tree is balanced, the distance to all 
the leaves is the same).
If $q \leq \ell+1$, then all $u,v \in \partial T$ are $R$ related. In
this case, it
is well known that that the topology of the tree may be recovered from
all minimal splits (this is the classical ``4 point method'', 
see e.g. \cite{ESSW1}).
We assume below that $q > \ell + 1$.  
Let $B_r(u) = \{v : d(v,u) = 2 r\}$.
Note that $B_r(u) \subset R(u)$, for all $u \in \partial T$ 
and $r \leq \ell+1$. 

\begin{Claim} \label{claim:splits}
$d$ satisfies 
\begin{itemize}
\item
d(u,v) = 0 if and only if $u = v$.
\item
For $1 \leq r \leq \ell$, 
$d(u,v) = 2 r$ if and only if $v \in R(u) \setminus B_{r-1}(u)$, and for all 
$\{w,w'\} \subset R(u) \setminus (B_{r-1}(u) \cup B_{r-1}(v))$,
the split $u v | w w'$ is a minimal split.
\end{itemize} 
\end{Claim}
\proofs
The first part is trivial.

For the second part note that
if $d(u,v) = 2r$, then $v \in R(u)$. Moreover, for all 
$w,w' \notin (B_{r-1}(u) \cup B_{r-1}(v))$, the intersection
$\pa(u,v) \cap \pa(w,w')$, is either empty or consists of a single
vertex. Therefore $u v | w w'$ is a minimal split.

If $d(u,v) < 2r$, then $v \notin R(u) \setminus B_{r-1}(u)$.

Suppose that $d(u,v) > 2r$ and $v \in  R(u)$. 
Since the tree is balanced, all the internal degrees are 
at least $3$ and $r+1 \leq \ell+1 < q$,
it follows that the sets $B_r(u) \setminus B_{r-1}(u)$, and  
$B_{r+1}(u) \setminus (B_r(u) \cup B_{r-1}(v))$ are not empty.
Let $u' \in B_r(u) \setminus B_{r-1}(u)$ and 
$v' \in  B_{r+1}(u) \setminus (B_r(u) \cup B_{r-1}(v))$.
Then $v',u' \in R(u) \setminus (B_{r-1}(u) \cup B_{r-1}(v))$ 
and $\path(u,u') \cap \path(v,v')$ 
is empty -- therefore $u v | u' v'$ is not a minimal split.
$\QED$

By Claim \ref{claim:splits}, from the minimal splits, 
we can recursively reconstruct for 
$r = 0, \ldots, \ell$ all pairs
 $u,v \in \partial T$ such that $d(u,v) = 2r$.
The second claim follows.
$\QED$

\prooft{of Theorem \ref{thm:4point}}
Note that by letting
\[
D(e) = \left\{ \begin{array}{ll} -\log(\theta(e)) & \mbox{ if } e \mbox{ is not adjacent to } \partial T, \\
                                 -\log(\theta(e) \eta(v)) & \mbox{ if } e = (u,v) \mbox{ and } v \in \partial T.
                \end{array} \right.
\]
the metric $D$ of Definition \ref{def:D} is of the form of the 
metric in Lemma \ref{lem:top}.

Moreover
\[
D_{\min} = 
\min \{D(e) : e \mbox{ not adjacent to } \partial T\} \geq 
\min_e -\log \theta(e) \geq -\log \theta_{\max}
> 1 - \theta_{\max}.
\]
Let $\eps' = - \log \theta_{\max}/4$ and $\eps = (1 - \theta_{\max})/4$.

We condition on the event that I. and II. of Lemma \ref{lem:d_D} hold with
$\eps$; so $2 D_{\min} \geq 8 \eps' = 8 \eps + 8(\eps' - \eps)$.
Thus for all $u$ and $v$ such that $d(u,v) \leq 2 \ell + 2$ 
it holds that $uR'v$.

Note that there exists a symmetric relation $\tilde{R}$ such that 
$R' \subset \tilde{R} \subset R$ and 
such that for
all $u$ and $v$ it is decidable in time polynomial in $k$ if they 
are $\tilde{R}$ related or not 
(to compute $\tilde{R}$ it suffices to check if 
$c(u,v) \geq \frac{3}{2} \theta_{\ast}$ within accuracy $\theta_{\ast}/4$). 

For a split $\Gamma = u_1 u_2 | u_3 u_4$, write $D^{\ast}(\Gamma)$ for 
$D^{\ast}(u_1,u_2) + D^{\ast}(u_3,u_4)$. 

Fix $u \in \partial T$ and 
$U = \{u_1,u_2,u_3,u_4\} \subset \tilde{R}(u)$. 
For all splits $\Gamma$ of
$U$, $|D^{\ast}(\Gamma) - D(\Gamma)| < 2 \eps$.

Let $\Gamma_1,\Gamma_2$ be splits
of $U$. We claim that $D(\Gamma_1) \leq D(\Gamma_2)$ if and only if
$D^{\ast}(\Gamma_1)  < D^{\ast}(\Gamma_2) + 4 \eps$ if and only if
$D^{\ast}(\Gamma_1)  < D^{\ast}(\Gamma_2) + 4 \eps'$.
Indeed, if $D^{\ast}(\Gamma_1)  < D^{\ast}(\Gamma_2) + 4 \eps'$,
then $D(\Gamma_1) < D(\Gamma_2) + 4 \eps' + 4 \eps 
<  D(\Gamma_2) + 2 D_{\min}$, 
and therefore $D(\Gamma_1) \leq D(\Gamma_2)$, by the first 
part of Lemma \ref{lem:top}. If on the other hand,
$D^{\ast}(\Gamma_1) \geq D^{\ast}(\Gamma_2) + 4 \eps$, then 
$D(\Gamma_1) > D(\Gamma_2)$, as needed.

Moreover, given that either
$D^{\ast}(\Gamma_1) \geq D^{\ast}(\Gamma_2) + 4 \eps'$ or
$D^{\ast}(\Gamma_1) < D^{\ast}(\Gamma_2)    + 4 \eps$, we may find
which of the two hold in time polynomial in $k$. 
Therefore, the minimal splits may be recovered in time polynomial in
$n$ and $k$.

We therefore conclude that conditioned on I. and II. of Lemma \ref{lem:d_D},
we may recover all the minimal splits of
$U$, for all $U \subset \tilde{R}(u)$ and all $u \in \partial T$ in
time polynomial in $k$ and $n$.

It now follows from the second part of Lemma \ref{lem:top} and 
from Lemma \ref{lem:d_D} that we 
may recover the $\ell$-topology of the tree with
error probability bounded by 
(\ref{eq:error_eps}):
\begin{eqnarray*}
n^2 \exp(-\frac{k}{8} \theta_{\ast}^4 \eps^2) &=&
n^2 \exp \left(-\frac{k}{8} 
    \left( \frac{\theta_{\min}^{2 \ell + 2} \eta_{\min}^2}{2} \right)^4
    \left( \frac{1 - \theta_{\max}}{4} \right)^2 \right)  \\ &=&
n^2 \exp \left( -\frac{k}{2048} \theta_{\min}^{8 \ell + 8} 
                \eta_{\min}^8 (1 - \theta_{\max})^2 \right),   
\end{eqnarray*}
as needed.

Finally note that
given the relation $\tilde{R}$ and all the minimal splits, the 
reconstruction procedure described in Lemma \ref{lem:top} is
computable in time polynomial in $n$. 
We conclude that 
the function
$\Phi$ is computable in time polynomial in $n$ and $k$. 

$\QED$


\section{Reconstruction of balanced trees} \label{sec:alg}

The proof of Theorem \ref{thm:CFN2} is similar to that of Theorem 
\ref{thm:CFN1}. The main difference is that instead of just
calculating correlations in order to recover $\ell$-topology, the $4$-point
method, i.e., Theorem \ref{thm:4point} is applied. The analysis
of the majority function in the more general setting, i.e., Theorem 
\ref{thm:maj_good} is also needed.

\prooft{of Theorem \ref{thm:CFN2}}
Let $b$ and $\theta_{\min}$ be such that $b \theta_{\min}^2 > 1$.
By Theorem \ref{thm:maj_good} there exist
$\ell,\alpha > 1$ and $\beta > 0$ be such that (\ref{eq:maj_prop}) holds.

To recover $d$, we will apply Theorem \ref{thm:4point} 
and Lemma \ref{cor:rec_maj} recursively in order to recover
$d^{\ast}_{i \ell}$, 
for $i = 0, \ldots, \lceil q/\ell \rceil$, where $q$ is the distance
from the root of the tree to the leaves. 

We note that the algorithms in Theorem \ref{thm:4point} and in 
Lemma \ref{cor:rec_maj} are polynomial time algorithms in $k$ and
$n$ -- since $k$ is polynomial in $n$, it follows 
that the running time of the reconstruction algorithm
below is polynomial in $n$.

Trivially, $d^{\ast}_0(v,u) = 2 \one_{v \neq u}$.
We show how given $d^{\ast}_{i \ell}$ and the samples 
$(\sigma^t_\partial)_{t=1}^k$, we can recover $d^{\ast}_{i \ell + \ell}$ with
error probability bounded by $n^2 \exp( - \tilde{c} k ) / b^{2 i
  \ell}$, where 
\begin{equation} \label{eq:tildec}
\tilde{c} = c^{\ast}\, \theta_{\min}^{8 \ell + 8} 
            \beta^8 (1 - \theta_{\max})^2,
\end{equation}
and $c^{\ast} \geq 1/2048$.

Let $\Psi_i : \{\pm 1\}^{\partial T} \to \{\pm 1\}^{L_{\partial - i
    \ell}}$ 
be the function defined in  
Lemma \ref{cor:rec_maj} given $d^{\ast}_{i \ell}$.
Then
\[
(\Psi_i(\sigma^t_{\partial T}))_{t=1}^k = 
(\sigma^t_v \tau^t_v : v \in L_{\partial - i \ell})_{t=1}^k,
\]
where $\tau^t_v$ are independent variables with $\E[\tau^t_v] \geq
\beta$. Moreover, $\tau^t_v$ are independent of 
$(\sigma^t_v : d(v,\partial T) \geq i \ell, 1 \leq t \leq k)$.

By Theorem \ref{thm:4point}, given 
$(\sigma^t_v \tau^t_v : v \in L_{\partial - i \ell})_{t=1}^k$, 
we may recover 
\[
d' : L_{\partial - i \ell} \times L_{\partial - i \ell} \to 
\{0,\ldots, 2 \ell + 2\},
\]
defined by $d'(u,v) = \min \{d(u,v), 2 \ell + 2\}$, 
with error probability bounded
by $n^2 \exp (-\tilde{c} k) / b^{2 i \ell}$.
As in Theorem \ref{thm:CFN1}, it easy to write $d^{\ast}_{i \ell + \ell}$ in
terms of $d^{\ast}_{i \ell}$ and $d'$.
  
Letting $A_i$ be the event of error in recovering
$d^{\ast}_{i \ell + \ell}$ given $d^{\ast}_{i \ell}$, 
and $\alpha = \sum_{i=0}^{\lceil q/l \rceil} \P[A_i]$, 
the total error probability is at most $\alpha$.

Now
\[
\alpha \leq \exp (-\tilde{c} k)
\left(n^2 + n^2/b^{2 \ell} + n^2/b^{4 \ell} + \cdots \right) \leq
2 n^2 \exp (-\tilde{c} k).
\] 
Defining ${c'}^{-1} = \tilde{c}$, and taking
\begin{equation} \label{eq:k_general}
k = \frac{\log(2 n^2) - \log \delta}{\tilde{c}} = 
c'(2 \log n + \log 2 - \log \delta)
\end{equation}
we obtain $\alpha \leq \delta$.
The statement of the theorem follows from (\ref{eq:k_general}) and
(\ref{eq:tildec}).
$\QED$

\prooft{of Theorem \ref{thm:CFN3}}
The proof is similar to that of Theorem \ref{thm:CFN2}.
Let $h^2 = (g^2 + b \theta_{\min}^2)/2$, so that $b \theta_{\min}^2 > h^2$. 
We choose $\ell,\alpha > h^{\ell}$ 
and $\beta > 0$ such that (\ref{eq:maj_prop}) holds.
The main difference from the proof of Theorem \ref{thm:CFN2} 
is that when we recover
$(\sigma^t_v \tau^t_v)_{v \in L_{\partial - i \ell}, 1 \leq t \leq k}$,
the $\tau^t_v$ are independent variables satisfying 
a weaker inequality, $\E[\tau^t_v] \geq \beta h^{i \ell}$.

Therefore
\[
\P[A_i] \leq \frac{n^2}{b^{2 i \ell}} \exp(-\tilde{c} h^{8 i \ell} k)
\leq \frac{n^2}{b^{2 i \ell}} \exp(-\tilde{c} h^{8 q} k),
\]
where $\tilde{c}$ is given in (\ref{eq:tildec}), and $q = \log_b n$.

If $k = c g^{-8q}$, then
\[
\sum \P[A_i] \leq 2n^2 \exp(-\tilde{c} c (h/g)^{8q}),
\] 
which is smaller than $\delta$ for all $n$, for sufficiently large 
value of $c$. 
$\QED$

{\bf Acknowledgments:} I wish to thank Mike Steel for proposing the
conjecture which motivated this work and for helpful comments on drafts
of this paper. Thanks to Noam Berger and Yuval Peres for helpful
discussions and to Lea Popovic for helpful comments on a draft of this
paper. Finally, many thanks to the anonymous referee for numerous of helpful 
suggestions and remarks.
 



\begin{thebibliography}{19}
\bibitem{Io} D. Ioffe (1996).  A note on the extremality of the disordered
state for the Ising model on the Bethe lattice.
{\em Lett.\  Math.\  Phys.} {\bf 37}, 137--143.
\end{thebibliography}
\end{document}